\documentclass[12pt]{article}

\usepackage{latexsym}
\usepackage{mathrsfs}
\usepackage{amsfonts,amsmath,amssymb}
\usepackage{indentfirst}
\usepackage{subeqnarray}
\usepackage{verbatim}
\usepackage[pdftex]{graphicx}
\usepackage{epstopdf}
\usepackage{subfigure}

\newcommand{\bx}{\mbox{\boldmath{$x$}}}

\newcommand{\bz}{\mbox{\boldmath{$z$}}}

\newcommand{\btau}{\mbox{\boldmath{$\tau$}}}
\newcommand{\bzero}{\mbox{\boldmath{$0$}}}

\newcommand{\bH}{\mbox{\boldmath{$H$}}}

\newcommand{\fb}{\mbox{\boldmath{$f$}}}

\newcommand{\bu}{\mbox{\boldmath{$u$}}}

\newcommand{\bv}{\mbox{\boldmath{$v$}}}

\newcommand{\bV}{\mbox{\boldmath{$V$}}}
\newcommand{\bw}{\mbox{\boldmath{$w$}}}

\newcommand{\real}{\mbox{$\mathbb{R}$}}

\newcommand{\bvarepsilon}{\mbox{\boldmath{$\varepsilon$}}}
\newcommand{\bsigma}{\mbox{\boldmath{$\sigma$}}}

\newcommand{\bnu}{\mbox{\boldmath{$\nu$}}}

\newcommand{\weak}{\rightharpoonup}

\newtheorem{theorem}{Theorem}[section]

\newtheorem{problem}[theorem]{Problem}

\newtheorem{proposition}[theorem]{Proposition}

\numberwithin{equation}{section}

\usepackage[pdftex,unicode,colorlinks,linkcolor=blue]{hyperref}

\begin{document}

\begin{center}
\Large\bf Variational-hemivariational inequalities:\\
A brief survey on mathematical theory and\\
numerical analysis
\end{center}

\medskip
\begin{center}
{\large\sc Weimin Han}\footnote{The work was partially supported by Simons Foundation Collaboration Grants (Grant No.\ 850737).}

Department of Mathematics

University of Iowa

Iowa City, IA 52242, USA

Email: {\tt weimin-han@uiowa.edu} 
\end{center}

\medskip
\begin{quote}
{\bf Abstract.}  Variational-hemivariational inequalities are an area full of interesting and 
challenging mathematical problems.  The area can be viewed as a natural extension of that of 
variational inequalities. Variational-hemivariational inequalities are valuable for application problems 
from physical sciences and engineering that involve non-smooth and even set-valued relations, monotone 
or non-monotone, among physical quantities. In the recent years, there has been substantial growth of 
research interest in modeling, well-posedness analysis, development of numerical methods and numerical 
algorithms of variational-hemivariational inequalities.  This survey paper is devoted to a brief account 
of well-posedness and numerical analysis results for variational-hemivariational
inequalities. The theoretical results are presented for a family of abstract stationary 
variational-hemivariational inequalities and the main idea is explained
for an accessible proof of existence and uniqueness. To better appreciate the distinguished feature of
variational-hemivariational inequalities, for comparison, three mechanical 
problems are introduced leading to a variational equation, a variational inequality, and a
variational-hemivariational inequality, respectively.  The paper also comments on mixed
variational-hemivariational inequalities, with examples from applications in fluid mechanics, and 
on results concerning the numerical solution of other types (nonstationary, history dependent)
of variational-hemivariational inequalities.
\end{quote}

\smallskip
{\bf Keywords:} Variational inequality, hemivariational inequality, variational-hemivariational inequality,
well-posedness, numerical solution, finite element method, discontinuous Galerkin method,
virtual element method, convergence, error estimate,
contact mechanics, Stokes hemivariational inequality, Navier-Stokes hemivariational inequality

{\bf AMS Mathematics Subject Classification (2020):} 49J40, 65N30, 74G15, 74G22, 74G30, 74M10, 74M15, 76D03, 76M10

\medskip
\section{Introduction}

It is generally agreed that the first variational inequality (VI) was formulated in A. Signorini's 
paper \cite{Si33} in 1933 for a study of a frictionless contact problem between a linearized elastic
body and a rigid foundation.  The problem is commonly known as the Signorini problem.  The existence 
and uniqueness of a solution to the Signorini problem was proved in early 1960s by G. Fishera
(cf.\ \cite{Fi64}).  This is considered as the beginning of the area of variational inequalities (VIs)  
(\cite{Ant83}).  In mid 1960s to early 1970s, foundations of basic mathematical theory of VIs were 
established in a series of papers, cf.\ \cite{Br72, HS66, LS67, St64}.  The monograph \cite{DL1976}
plays an important role in popularizing the area of VIs as it is shown in the book that many complicated 
application problems in mechanics and physics can be modeled and studied as VIs. Since there are no 
analytic solution formulas for VIs arising in applications, one relies on numerical methods to solve VIs.  
Early comprehensive references on numerical methods for solving VIs are \cite{Gl1984, GLT1981, HHNL1988}. 
Mechanics is a rich source of VIs, and there is a large number of references devoted to this topic, 
cf.\ \cite{HR2013, Te1985} on elasto-plasticity, 
\cite{CHR2023, EJK2005, HS2002, HHN96, KO1988, SHS2006, SM2009, SM2012} 
on contact mechanics.  Nowadays, active research persists in the area of VIs due to emerging new
applications and the need of developing more efficient and effective numerical methods and algorithms 
to solve VIs (e.g., \cite{CHY23, GJKR2022, JA2023, UL2011, Yo21}).

VIs are featured by the presence of nonsmooth terms with a convex structure in their mathematical 
formulations.  For applications, the nonsmooth convex structure often comes from a nonsmooth 
monotone relation among physical quantities of interest.  For applications involving nonsmooth
nonmonotone relations for physical quantities, hemivariational inequalities (HVIs) arise.  
P. D. Panagiotopoulos started the area of HVIs in early 1980s (\cite{Pa83}).  VIs can be viewed as
a special or degenerate case of HVIs in the sense that when the nonsmooth relations among physical 
quantities happen to be monotone, a HVI is reduced to a VI.  We may consider the more general
variational-hemivariational inequality (VHI) which contains nonsmooth terms of both kinds, those
with a convex structure and those with a nonconvex structure.  In the literature, the term HVI is 
often used to refer to a VHI also.  In this paper, we use the terms VHI and HVI interchangeably.

Early comprehensive references in the area of VHIs include \cite{MP1999, NP1995, Pa1993}
on modeling, mathematical analysis and applications, and \cite{HMP1999} for the finite element method 
to solve VHIs.  In the last three decades, the area of VHIs has attracted the attention of ever more 
researchers, and recent comprehensive coverage of mathematical theory and applications can be found in
\cite{CL2021, CLM2007, GM2003, GMDR2003, MOS2013, SM2025}.  In these books and most math journal papers
on VHIs, abstract surjectivity results on pseudomonotone operators are needed in proving the existence of 
a solution. Such an approach has its own merits.  However, the requirement of the knowledge on abstract 
mathematical theory of pseudomonotone operators seems to be a hurdle to popularize the area of VHIs
in the research communities of applied and computational mathematicians and engineers. An effort was 
made in \cite{Han20, Han21} in developing an alternative accessible approach for the mathematical theory 
of VHIs that does not rely on the abstract mathematical theory of pseudomonotone operators. We will 
describe the main idea of this alternative approach for studies of stationary VHIs in 
Subsection \ref{subsec:abs}.  We also note that the accessible approach is extended in \cite{HM22a, HM22b} 
for well-posedness of mixed VHIs.  A comprehensive reference in these regards is the book \cite{Han2024}.

Numerical methods are needed to solve VHIs. The finite element method and a variety of solution 
algorithms are discussed in \cite{HMP1999} to solve HVIs. An optimal order error estimate is first 
presented in \cite{HMS14} for the linear finite element solutions of a VHI.  This is followed by 
a series of papers on further analysis of the finite element method to solve VHIs, e.g., 
\cite{Han18, HSB17, HSD18, HZ19}.  The survey papers \cite{HS19AN} and \cite{HFWH25} provide accounts 
of recent advances on numerical analysis of VHIs.  Besides the finite element method, 
other numerical methods such as the discontinuous Galerkin method and the virtual element method, 
have been also been applied to solve VHIs, cf.\ discussions in Section \ref{sec:other}.
For the numerical solution of inequality problems (VIs, VHIs) of second-order, an optimal order error 
estimate can be achieved only for the linear element (linear finite element, linear virtual element) solutions
(order one in the $H^1$ norm).  Moreover, due to the limited regularity properties for solutions of 
inequality problems, low order elements are preferred to solve the inequalities.  

The paper provides a summary account on the mathematical theory and numerical solution of VHIs and we
focus on the stationary/time-independent problems.  The main goal is to present the reader a relatively
complete picture on the current status of the research on VHIs, especially regarding the numerical 
analysis of VHIs.  In Section \ref{sec:pre}, we review the notions of generalized subdifferentials and 
generalized subgradients, and their properties, as these are fundamental in the study of VHIs.  
In Section \ref{sec:sample}, we introduce three typical sample problems in mechanics in the 
order of increasing complexity. The first example is a variational equation (VE) for a standard 
boundary value problem. The other two examples are from contact mechanics, in the form of a VI
and a VHI, respectively. It is hoped that through the presentation of the three examples, one
can observe clearly what kind of features in a mechanical problem lead to a VI or a VHI.
In Section \ref{sec:abs}, we consider an abstract stationary VHI.  The first part of the section 
is to explain the main idea in the well-posedness analysis without the commonly used surjectivity 
results of pseudomonotone operators.  The second part of the section summarizes main numerical
analysis results on the abstract stationary VHI. The results presented in this section are general,
and they can be applied to concrete VHIs, and to VIs and VEs which are special cases of VHIs. 
In Section \ref{sec:contact}, the theory presented in Section \ref{sec:abs} is applied to
the VHI contact problem introduced in Section \ref{sec:sample}, and we provide a well-posedness 
result and an optimal order error estimate for its numerical solution using the linear finite 
element method  under certain solution regularity assumptions.  In Section \ref{sec:fluid}, 
we present sample results on Stokes and Navier-Stokes HVIs in fluid mechanics, as examples 
of mixed VHIs.  In Section \ref{sec:other}, we briefly comment on other numerical methods to 
solve VHIs and studies of time-dependent VHIs.

\section{Generalized directional derivative, subdifferential and subgradient}\label{sec:pre}

In the study of nonsmooth problems, one fundamental issue is how to extend the concept of differentiability 
for functions that are not differentiable in the classical sense. For VIs, notions of subdifferental and
subgradients serve such a purpose, and one can consult \cite{ET1976, ZeIII} and many other excellent 
books on convex analysis for these notions.  For VHIs, we need the notions of the generalized directional 
derivative and generalized subdifferential/subgradient for locally Lipschitz continuous functions
introduced by F. H. Clarke (\cite{Cl75, Cl1983}). 

Let $V$ be a Banach space, and let $U$ be an open subset in $V$. Often, we can simply take $U=V$.  
Let $\Psi\colon U\to \real$ be a locally Lipschitz continuous function. Given $u\in U$ and $v\in V$,
the classical directional directive $\Psi^\prime(u;v)$ may not exist.  The idea is to consider the 
ratio $\left(\Psi(w+\lambda v) -\Psi(w)\right)/\lambda$ for $w$ close to $u$ and $\lambda>0$ close to 0.
If the classical directional directive $\Psi^\prime(u;v)$ does not exist, then we cannot take the limit 
as $w\to u$ and $\lambda\downarrow 0$ on the ratio $\left(\Psi(w+\lambda v) -\Psi(w)\right)/\lambda$.
Nevertheless, since $\Psi$ is locally Lipschitz continuous at $u$, the upper limit of the ratio always 
exists as a real number.  Thus, the generalized (or Clarke) directional derivative of $\Psi$ at $u\in U$
in the direction $v \in V$ is defined by
\[\Psi^0(u; v) := \limsup_{w\to u,\,\lambda \downarrow 0}\frac{\Psi(w+\lambda v) -\Psi(w)}{\lambda}.\]
The next step is to define the generalized (or Clarke) subdifferential of $\Psi$ at $u\in U$ by
\begin{equation}
\partial\Psi(u):=\left\{u^*\in V^*\mid\Psi^{0}(u;v)\ge\langle u^*,v\rangle\ \forall\,v\in V\right\}. 
\label{gen:sub}
\end{equation}
It can be shown that $\partial \Psi(u)$ is nonempty, convex, and weak\-ly$^{\,*}$ compact in $V^*$.
Any element in the set $\partial \Psi(u)$ is called a generalized (or Clarke) subgradient of $\Psi$ at 
$u\in U$.  

Furthermore, it can be proved that in case $\Psi \colon U \to \real$ happens to be convex, 
then the generalized subdifferential $\partial \Psi(u)$ at any $u \in U$ coincides with the 
convex subdifferential $\partial \Psi(u)$.  Because of this property, it is perfectly reasonable
to use the same symbol $\partial$ for both the generalized subdifferential of locally Lipschitz 
continuous functions and the convex subdifferential of convex functions.

Properties of the generalized directional derivative and the generalized subdifferential can be found in 
several books, e.g., \cite{Cl1983, CLSW1998} or \cite[Section 3.2]{MOS2013}.  In the following, 
we mention some basic properties to help those readers without prior exposures to the notions of the 
generalized directional derivative and the generalized subdifferential for a better understanding.

The generalized directional derivative is non-negatively homogeneous and sub-additive with respect to 
the direction variable:
\begin{align*}
& \Psi^0(u;\lambda\,v)= \lambda\, \Psi^0(u; v)\quad\forall\,\lambda\ge 0, \,u\in U,\,v\in V,\\
& \Psi^0 (u; v_1 + v_2) \le \Psi^0(u; v_1)+\Psi^0(u; v_2)\quad\forall\,u\in U,\,v_1,v_2\in V.
\end{align*}
The generalized subdifferential is defined by the formula \eqref{gen:sub} through the use of the 
generalized directional derivative.  Conversely, knowing the generalized subdifferential, 
we can compute the generalized directional derivative by the formula
\begin{equation}
\Psi^0(u;v)=\max\left\{\langle u^*,v\rangle\mid u^*\in\partial\Psi(u)\right\}\quad\forall\,u\in U,\,v\in V.
\label{gen:dir}
\end{equation}

In the study of VHIs and their numerical approximations, limiting properties of 
the generalized subdifferential and the generalized directional derivative are useful 
and we list two such properties below:\\  
If $u_n\to u$ in $V$, $u_n,u\in U$, and $v_n \to v$ in $V$, then 
\[ \limsup_{n\to\infty} \Psi^0(u_n;v_n) \le \Psi^0(u; v). \]
If $u_n \to u$ in $V$, $u_n,u\in U$, $u^*_n \in \partial \Psi(u_n)$,
and $u^*_n \to u^*$ weakly$^{*}$ in $V^*$, then $u^*\in\partial\Psi(u)$.

Next, let $\Psi,\Psi_1,\Psi_2 \colon U\to \real$ be locally Lipschitz functions. Then, we have 
the scalar multiplication rule
\[ \partial(\lambda\,\Psi)(u)=\lambda\,\partial \Psi(u)\quad \forall\,\lambda \in \real, \ u\in U,\]
and the summation rule
\begin{equation}\label{B.PART11}
\partial (\Psi_1 + \Psi_2) (u) \subset \partial \Psi_1 (u) + \partial \Psi_2 (u)\quad\forall\,u\in U,
\end{equation}
or equivalently,
\begin{equation}\label{B.PART12}
(\Psi_1 + \Psi_2)^0(u; v) \le \Psi_1^0(u; v) + \Psi_2^0(u; v)\quad\forall\,u\in U,\, v\in V.
\end{equation}
Moreover, \eqref{B.PART11} and \eqref{B.PART12} hold with equalities if $\Psi_1$ and $\Psi_2$ are 
regular at $u$.  The regularity of the Lipschitz continuous function $\Psi\colon U\to\mathbb{R}$
at $u\in U$ means that the directional derivative $\Psi^\prime(u;v)$ exists and 
\begin{equation}
\Psi^\prime(u;v) = \Psi^0(u;v)\quad\forall\,v\in V.
\label{eq:regular}
\end{equation}
It is known that a function is regular at any point where the function is continuously differentiable.
Also, a convex function is regular in the interior of its effective domain.

In the study of VHIs, we will assume that there exists a constant $\alpha_\Psi\ge 0$ such that
\begin{equation}
\Psi^0(v_1;v_2-v_1) + \Psi^0(v_2;v_1-v_2) \le \alpha_\Psi \|v_1-v_2\|_V^2 \quad\forall\,v_1,v_2\in U.
\label{2.12a}
\end{equation}
This condition characterizes the degree of non-convexity of the functional $\Psi$: the smaller the 
constant $\alpha_\Psi\ge 0$, the weaker the non-convexity of $\Psi$.
For a convex functional $\Psi$, \eqref{2.12a} holds with $\alpha_\Psi=0$.  The condition \eqref{2.12a}
is sometimes expressed equivalently as a condition on the generalized subdifferential: 
\begin{equation}
\langle v^*_1-v^*_2,v_1-v_2\rangle \ge -\alpha_\Psi \|v_1-v_2\|_V^2 \quad\forall\,v_i\in U,
\,v^*_i\in \partial\Psi(v_i),\, i=1,2.
\label{2.12b}
\end{equation}
In the literature, \eqref{2.12b} is usually called a relaxed monotonicity condition. The inequality 
\eqref{2.12b} with $\alpha_\Psi=0$ is the monotonicity of $\partial\Psi$ for a convex functional $\Psi$.
A short-hand expression of the condition \eqref{2.12b} is 
\begin{equation}
\langle \partial\Psi(v_1)-\partial\Psi(v_2),v_1-v_2\rangle\ge -\alpha_\Psi\|v_1-v_2\|_V^2 \quad\forall\,v_1,v_2\in U.
\label{2.12c}
\end{equation}
It can be proved (e.g., \cite[p.\ 26]{Han2024}) that the condition \eqref{2.12a} holds if and only if 
the functional $v\mapsto \Psi(v)+(\alpha_\Psi/2)\,\|v\|_V^2$ is convex on $U$.  This result provides 
a simple way to verify the condition \eqref{2.12a}.

\smallskip
In virtually all the applications in mechanics, the locally Lipschitz continuous functional $\Psi$ takes
the form of an integral of a locally Lipschitz continuous function $\psi$ of a real variable or of 
several real variables. For a locally Lipschitz continuous function defined over a finite dimensional set, 
there is a useful formula to compute the generalized subdifferential 
(cf.\ \cite[Theorem 10.7]{Cl2013}, \cite[Prop.\ 3.34]{MOS2013}).

\begin{proposition} \label{prop:finite-dim}
Assume $\Omega\subset \mathbb{R}^d$ is open, $\psi\colon \Omega \to\mathbb{R}$ is locally Lipschitz 
continuous near $\bx\in \Omega$, $D\subset\mathbb{R}^d$ with $|D|=0$, and $D_\psi\subset\mathbb{R}^d$ with
$|D_\psi|=0$ such that $\psi$ is Fr\'{e}chet differentiable on $\Omega\backslash D_\psi$.  Then,
\[ \partial\psi(\bx)={\rm conv}\left\{\lim \psi^\prime(\bx_k)\mid \bx_k\to\bx,\
 \bx_k\not\in D\cup D_\psi\right\}. \]
\end{proposition}

Next, we show some examples to compute the generalized subdifferential of locally Lipschitz 
continuous functions by applying Proposition \ref{prop:finite-dim} to conclude this section. 

For the function $\psi_1(x)=-|x|$, $x\in\mathbb{R}$, its generalized subdifferential coincides with its 
classical derivative $-1$ for $x>0$, and $1$ for $x<0$.  At $x=0$, $\partial\psi_1(0)$ equals the 
convex hull of the two derivative limiting points $-1$ and $1$.  Thus,
\[ \partial\psi_1(x)=\left\{\begin{array}{ll} 1 & {\rm if}\ x<0,\\ \lbrack -1,1 \rbrack & {\rm if}\ x=0,\\
-1 & {\rm if}\ x>0.\end{array} \right. \]
For the generalized directional directive $\psi_1(x;y)$, $x,y\in\mathbb{R}$, we can use the 
property \eqref{gen:dir} to find that
\[ \psi^0_1(x;y)=\left\{\begin{array}{ll} y & {\rm if}\ x<0,\\ |y| & {\rm if}\ x=0,\\
-y & {\rm if}\ x>0.\end{array} \right. \]

Similarly, for $\psi_2(x)=|x|$, $x\in\mathbb{R}$, its generalized subdifferential is
\[ \partial\psi_2(x)=\left\{\begin{array}{ll} -1&{\rm if}\ x<0,\\ \lbrack -1,1\rbrack & {\rm if}\ x=0,\\
1& {\rm if}\ x>0. \end{array}\right. \]
Note that $\psi_2$ is a convex function, and $\partial\psi_2$ is also the convex subdifferential of $\psi_2$.  Moreover, the generalized directional directive is, for any $y\in\mathbb{R}$,
\[ \psi^0_2(x;y)=\left\{\begin{array}{ll} -y & {\rm if}\ x<0,\\ |y| & {\rm if}\ x=0,\\
y & {\rm if}\ x>0.\end{array} \right. \]

As a more complicated example, consider
\[ 
\psi_3(x)=\left\{\begin{array}{ll} \sin x & {\rm if}\ x< 0,\\[0.1mm]
-x^2  & {\rm if}\ 0\le x \le 1,\\[0.1mm] x^3-3\,x^2+3\,x-2 & {\rm if}\ x>1.
\end{array}\right.
\]
For its generalized subdifferential, we can readily write down the formula
\[
\partial \psi_3(x)=\left\{\begin{array}{ll} \cos x& {\rm if}\ x< 0,\\[0.1mm]
[0,1]& {\rm if}\ x=0,\\[0.1mm]
-2\,x & {\rm if}\ 0<x<1,\\[0.1mm]
[-2,0] & {\rm if}\ x=1,\\[0.1mm]
3\left(x-1\right)^2 & {\rm if}\ x>1.
\end{array}\right.
\]
For the directional derivative at points where the classical derivative does not exist, we find that
\begin{align*}
\psi^0_3(0;y) & = \max\{0,y\},\\
\psi^0_3(1;y) & = \max\{0,-2\,y\}
\end{align*}
for any $y\in\mathbb{R}$.

\section{Three representative problems in linearized elasticity}\label{sec:sample}

We will see an example of a VHI in Subsection \ref{sec:VHI}, for a contact problem.  For comparison,
under similar mechanical settings with simpler boundary conditions, we will provide an example of a 
variational equation in Subsection \ref{sec:VE}, and an example of a VI in Subsection \ref{sec:VI}. 
We will first introduce the notation in Subsection \ref{subsec:notation}.

\subsection{Notation}\label{subsec:notation}

Consider mathematical models describing the equilibrium state of an elastic body subject to the action 
of external forces and constraints on the boundary.  Let the reference configuration of the body 
be the closure of an open, bounded and connected set $\Omega$ in $\real^d$. For applications, 
the dimension $d=2$ or $3$. Assume $\Omega$ has a Lipschitz boundary $\Gamma=\partial\Omega$. 
Then, the unit outward normal vector $\bnu$ exists a.e.\ on $\Gamma$.  We use boldface letters for 
vectors and tensors. A typical point in $\mathbb{R}^d$ is denoted by $\bx=(x_i)$. The range of indices 
$i$, $j$, $k$, $l$ is between $1$ and $d$.  We adopt the summation convention over a repeated index.  
For convenience, for a subset $\Delta$ of $\Omega$ or that of $\Gamma$, and a function $g$ defined on
$\Delta$, we use $I_\Delta(g)$ to denote the integral of the function $g$ over $\Delta$.

We denote by $\mathbb{S}^d$ the space of real symmetric matrices of order $d$. 
Over $\mathbb{R}^d$ and $\mathbb{S}^{d}$, we use the canonical inner products and norms:
\begin{align*}
& \bu\cdot\bv=u_i v_i,\quad|\bv|=(\bv\cdot\bv)^{1/2}\quad\forall\,\bu=(u_i),\bv=(v_i)\in\mathbb{R}^d,\\[0.1mm]
&\bsigma:\btau =\sigma_{ij}\tau_{ij},\quad |\btau|=(\btau:\btau)^{1/2}\quad\forall\,
\bsigma=(\sigma_{ij}),\btau=(\tau_{ij}) \in\mathbb{S}^{d}.  
\end{align*}

The primary unknown of the mechanical problems in this section is the displacement of the elastic body,
$\bu\colon\overline{\Omega}\to \mathbb{R}^d$.  The linearized strain tensor is
\[ \bvarepsilon(\bu)=\frac 12\left(\nabla \bu+(\nabla \bu)^T\right), \]
which is an $\mathbb{S}^d$-valued field in $\Omega$.  In componentwise form,
\begin{equation*}
\varepsilon_{ij}(\bu)=(\bvarepsilon(\bu))_{ij} =\frac 12\, (u_{i,j} + u_{j,i}),\quad 1\le i,j\le d,
\end{equation*}
where $u_{i,j}=\partial u_i/\partial x_j$. In mechanical problems, another important quantity is the 
stress tensor $\bsigma\colon \Omega\to \mathbb{S}^d$, which is also an $\mathbb{S}^d$-valued field in $\Omega$.

We will use Sobolev and Lebesgue spaces on $\Omega$, $\Gamma$, or their subsets, such as $L^2(\Omega;\real^d)$,
$L^2(\Gamma_N;\real^d)$, $L^2(\Gamma_C;\real^d)$, $H^1(\Omega;\real^d)$, and $H^1(\Omega;\mathbb{S}^d)$,
endowed with their canonical inner products and associated norms. Here, $\Gamma_N$ and $\Gamma_C$ are
measurable subsets of $\Gamma$. For a function $\bv\in H^1(\Omega;\real^d)$ we write $\bv$ for its trace 
$\gamma\bv\in L^2(\Gamma;\real^d)$ on $\Gamma$.  A standard reference on Sobolev spaces is \cite{AF2003}.  
One may also consult many other books on Sobolev spaces or PDEs, e.g., \cite{Bre2011, Evans2010}.

The space for the stress field is
\begin{equation}
\mathbb{Q}:=L^2(\Omega;\mathbb{S}^d)=\left\{\bsigma=(\sigma_{ij})\mid\sigma_{ij}=\sigma_{ji}\in L^2(\Omega),
\ 1\le i,j\le d\right\}.  \label{SpQ}
\end{equation}
This is a real Hilbert space endowed with the inner product
\[  (\bsigma,\btau )_\mathbb{Q}=\int_{\Omega}\bsigma:\btau\,dx,\quad \bsigma,\btau\in \mathbb{Q}.\]
The corresponding norm is denoted by $\|\cdot\|_\mathbb{Q}$. 

Let $\Gamma_D$ be a non-trivial measurable subset of $\Gamma$.  We will specify a homogeneous 
displacement boundary condition on $\Gamma_D$, and seek the unknown displacement field in the space
\begin{equation}
 \bV:=\left\{\bv\in H^1(\Omega;\mathbb{R}^d) \mid \bv=\bzero\ \mbox{on}\ \Gamma_D\right\}
\label{SpV}
\end{equation}
or its subspace.  Since $|\Gamma_D|>0$, Korn's inequality asserts that there is
a constant $c>0$, depending on $\Omega$ and $\Gamma_D$, such that
\begin{equation}
\|\bv\|_{H^1(\Omega;\mathbb{R}^d)}\le c\,\|\bvarepsilon(\bv)\|_\mathbb{Q}\quad\forall\,\bv\in \bV.
\label{Korn}
\end{equation}
A proof of the Korn inequality can be found in numerous publications, 
e.g.\ \cite[p.\ 79]{NH1981}.  As a result, $\bV$ is a real Hilbert space under the inner product
\begin{equation}
(\bu,\bv)_{\boldsymbol V}=(\bvarepsilon(\bu),\bvarepsilon(\bv))_\mathbb{Q}. \label{3.3}
\end{equation}
The induced norm is
\begin{equation}\label{3.3n}
\|\bv\|_{\boldsymbol V}=\|\bvarepsilon(\bv)\|_\mathbb{Q}.
\end{equation}
It follows from Korn's inequality \eqref{Korn} that $\|\cdot\|_{H^1(\Omega;\mathbb{R}^d)}$ and 
$\|\cdot\|_{\boldsymbol V}$ are equivalent norms on $\bV$. We will use $\|\cdot\|_{\boldsymbol V}$ 
as the norm on $\bV$ or its subspace.  Denote by $\bV^*$ the dual of the space $\bV$ and by 
$\langle\cdot,\cdot\rangle$ the corresponding duality pairing. For any element $\bv\in \bV$, define 
its normal and tangential components on $\Gamma$ by $v_\nu=\bv\cdot\bnu$ and $\bv_\tau=\bv-v_\nu\bnu$,
respectively. Similarly, for a function $\bsigma:\overline\Omega\to\mathbb{S}^d$, its normal and 
tangential components on $\Gamma$ are $\sigma_{\nu}=(\bsigma\bnu)\cdot\bnu$ and
$\bsigma_{\tau} =\bsigma\bnu - \sigma_{\nu}\bnu$, respectively.

For a differentiable field $\bsigma\colon \Omega\to\mathbb{S}^d$, its divergence is a vector-valued function
${\rm div}\,{\bsigma}\colon \overline{\Omega}\to\real$ with components
\[ ({\rm div}\,{\bsigma})_i=\sigma_{ij,j},\quad 1\le i\le d,\]
where the summation convention over the repeated index $j$ is applied.  
For $\bsigma\in H^1(\Omega;\mathbb{S}^d)$ and $\bv\in H^1(\Omega;\real^d)$, 
we have Green's formula (integration-by-parts formula)
\[  \int_\Omega\,\bsigma:\bvarepsilon(\bv)\,dx+\int_\Omega\,{\rm div}\,\bsigma\cdot\bv\,dx
= \int_\Gamma\bsigma\bnu \cdot\bv\,ds, \] 
which is applied repeatedly in derivations of weak formulations.

\subsection{A variational equation problem in linearized elasticity}\label{sec:VE}

The first mechanical problem leads to a variational equation.  In this subsection only, the boundary 
$\Gamma=\Gamma_D\cup\Gamma_N$ is split into two parts $\Gamma_D$ and $\Gamma_N$, with $|\Gamma_D|>0$. 
If $\Gamma_N=\emptyset$, then $\Gamma_D=\Gamma$ is the entire boundary.  
The pointwise formulation of the problem is 
\begin{align}
-\operatorname*{div}{\bsigma} & =\fb_0\quad\mathrm{in\ }\Omega, \label{ve1}\\
\bsigma & ={\cal E}\bvarepsilon(\bu)\quad\mathrm{in\ }\Omega, \label{ve2}\\
\bvarepsilon(\bu)&=\frac{1}{2}\left[\nabla\bu+(\nabla\bu)^T\right]\quad\mathrm{in\ }\Omega, \label{ve3}\\
\bu & =\mathbf{0}\quad\text{on}\ \Gamma_D, \label{ve4}\\
\bsigma{\bnu} & =\fb_2\quad\text{on}\ \Gamma_N. \label{ve5}
\end{align}
We comment that \eqref{ve1} is the equilibrium equation, $\fb_0$ being the density function of the 
body force; \eqref{ve2} is the elastic constitutive law, ${\cal E}$ being the elasticity tensor;
\eqref{ve3} is the defining relation for the linearized strain tensor; \eqref{ve4} represents the 
homogeneous boundary condition on $\Gamma_D$; and \eqref{ve5} describes the traction boundary condition 
on $\Gamma_N$, $\fb_2$ being the density function of the traction force.

In the general case, the elasticity operator ${\cal E}\colon \Omega\times \mathbb{S}^d \to \mathbb{S}^d$ 
in the constitutive law \eqref{ve2} is allowed to depend on the spatial location.  For homogeneous 
materials, ${\cal E}\colon \mathbb{S}^d \to \mathbb{S}^d$ does not depend on the spatial variable.
We assume the following properties:
\begin{equation}
\left\{\begin{array}{ll}
{\rm (a)\  There\ exists\ a\ constant}\ L_{\cal E}>0\ {\rm such\ that\ a.e.\ in}\ \Omega,\\
{}\qquad |{\cal E}(\cdot,\bvarepsilon_1)-{\cal E}(\cdot,\bvarepsilon_2)|
\le L_{\cal E} |\bvarepsilon_1-\bvarepsilon_2| \quad\forall\, \bvarepsilon_1,\bvarepsilon_2\in \mathbb{S}^d; \\ [1mm]
{\rm (b)\  there\ exists\ a\ constant}\ m_{\cal E}>0\ {\rm such\ that\ a.e.\ in}\ \Omega,\\
{}\qquad ({\cal E}(\cdot,\bvarepsilon_1)-{\cal E}(\cdot,\bvarepsilon_2)):
(\bvarepsilon_1-\bvarepsilon_2)\ge m_{\cal E} |\bvarepsilon_1-\bvarepsilon_2|^2\\
{}\qquad\qquad \forall\, \bvarepsilon_1,\bvarepsilon_2 \in \mathbb{S}^d; \\ [1mm]
{\rm (c) } \ {\cal E}(\cdot,\bvarepsilon)\ {\rm is\ measurable\ on\ }\Omega
\ {\rm for\ all \ }\bvarepsilon\in \mathbb{S}^d;  \\ [1mm]
{\rm (d)}\ {\cal E}(\cdot,\bzero)=\bzero\  {\rm a.e.\ in\ } \Omega.
\end{array}\right.
\label{Ass:E}
\end{equation}
For the force densities, we assume
\begin{equation}
\fb_0\in L^2(\Omega;\mathbb{R}^d), \quad \fb_2\in L^2(\Gamma_N;\mathbb{R}^d),
\label{Ass:f}
\end{equation}
and define and element $\fb\in\bV^*$ by
\begin{equation}
\langle \fb,\bv\rangle =\int_{\Omega}\fb_0\cdot\bv\,dx+\int_{\Gamma_N}\fb_2\cdot\bv\,ds\quad\forall\,\bv\in \bV.
\label{ell}
\end{equation} 

Through a standard procedure, one can derive the following weak formulation of the problem 
defined by \eqref{ve1}--\eqref{ve5}; for such a derivation, cf.\ e.g., \cite[Section 4.2]{Han2024}.

\begin{problem}\label{prob:ve}
{\it Find a displacement field $\bu\in \bV$ such that}
\begin{equation}
({\cal E}(\bvarepsilon({\bu})),\bvarepsilon(\bv))_\mathbb{Q}=\langle \fb,\bv\rangle\quad\forall\,\bv\in \bV,
\label{ve7}
\end{equation}
where $\fb\in\bV^*$ is defined by \eqref{ell}.
\end{problem}

Well-posedness of Problem \ref{prob:ve} can be shown through an application of the well-known 
Lax-Milgram lemma; for detail, one may consult \cite[Theorem 4.1]{Han2024}.

\subsection{A variational inequality in contact mechanics}\label{sec:VI}

In the study of contact problems, we assume the boundary of $\Gamma$ of the domain $\Omega$ is decomposed 
into three non-overlapping measurable parts: $\Gamma=\Gamma_D\cup\Gamma_N\cup\Gamma_C$, with 
$|\Gamma_D|>0$ and $|\Gamma_C|>0$.  We will specify a displacement boundary condition on $\Gamma_D$, 
a traction boundary condition on $\Gamma_N$, and contact boundary conditions on $\Gamma_C$. A variety 
of mathematical models of contact problems can be found in many publications, cf.\ e.g., the 
comprehensive references \cite{EJK2005, HS2002, KO1988, MOS2013, SM2009, SM2012, Wr2006}.

For the sample contact problem considered in this subsection, the pointwise relations 
\eqref{ve1}--\eqref{ve5} are supplemented by a bilateral contact condition with a Tresca's friction law on $\Gamma_C$:
\begin{equation}
\begin{array}
[c]{l}%
u_{\nu}=0,\\
|\bsigma_\tau|\le f_b\quad\mathrm{and}\\
{}\quad|\bsigma_\tau|<f_b\ \Rightarrow\ \bu_\tau=\bzero,\\
{}\quad|\bsigma_\tau|=f_b\ \Rightarrow\ \bu_\tau
=-\lambda\,\bsigma_\tau\ \mathrm{for\ some\ }\lambda\ge0.
\end{array}
\label{VI:mod1}
\end{equation}
Here $f_b\ge 0$ is a constant bound of the magnitude of the friction force.  

A study of this contact problem can be found in several references, e.g., \cite[Subsection 4.4.1]{Han2024}.
The tangential contact conditions in \eqref{VI:mod1} can be equivalently expressed as 
\[ |\bsigma_\tau|\le f_b\quad {\rm and}\quad \bsigma_\tau{\cdot}\bu_\tau+f_b |\bu_\tau|=0, \]
or as
\[ |\bsigma_\tau|\le f_b,\quad -\bsigma_\tau=f_b\,\frac{\bu_\tau}{|\bu_\tau|}\ {\rm if}\ \bu_\tau\not=\bzero,\]
or as
\[ -\bsigma_\tau\in\partial\phi(\bu_\tau),\quad \phi(\bz)=f_b|\bz|\ {\rm for}\ \bz\in\mathbb{R}^d.\]
Here, $\phi$ is a convex function and $\partial\phi$ stands for the (convex) subdifferential of $\phi$.

Define a subspace of the space $\bV$:
\begin{equation}
\bV_1=\left\{\bv\in\bV\mid v_\nu=0\ {\rm on}\ \Gamma_C\right\}
\label{SpV1}
\end{equation}
and use the norm $\|\cdot\|_{\boldsymbol V}$ over the subspace $\bV_1$. 
The weak formulation of the contact problem is the following VI (cf.\ \cite[Problem 4.6]{Han2024}).  

\begin{problem}\label{prob:VI1}
Find $\bu\in \bV_1$ such that
\begin{equation}
({\cal E}\bvarepsilon(\bu),\bvarepsilon(\bv-\bu))_{\mathbb{Q}}+I_{\Gamma_C}(f_b|\bv_\tau|)
-I_{\Gamma_C}(f_b |\bu_\tau|) \ge \langle \fb,\bv-\bu\rangle\quad\forall\,\bv\in \bV_1,
\label{1.50b}
\end{equation}
where $\fb\in\bV^*$ is defined by \eqref{ell}.
\end{problem}

Well-posedness of this problem is the content of \cite[Theorem 4.7]{Han2024}.  
See also Section \ref{sec:contact}.

\subsection{A variational-hemivariational inequality in contact mechanics}\label{sec:VHI}

In this subsection, we consider another contact problem whose mathematical model is a VHI.
As in the previous subsection, $\Gamma=\Gamma_D\cup\Gamma_N\cup\Gamma_C$, with $|\Gamma_D|>0$ and 
$|\Gamma_C|>0$. The pointwise relations \eqref{ve1}--\eqref{ve5} are supplemented by the normal 
compliance contact condition (\cite[Section 6.3]{MOS2013}) with Tresca's friction law on the contact 
boudnary $\Gamma_C$:
\begin{equation}
\begin{array}
[c]{l}
-\sigma_\nu\in\partial \psi_\nu(u_\nu),\\
|\bsigma_\tau|\le f_b\quad\mathrm{and}\\
{}\quad|\bsigma_\tau|<f_b\ \Rightarrow\ \bu_\tau=\bzero,\\
{}\quad|\bsigma_\tau|=f_b\ \Rightarrow\ \bu_\tau
=-\lambda\,\bsigma_\tau\ \mathrm{for\ some\ }\lambda\ge0.
\end{array}
\label{cont6}
\end{equation}
Here, the function $\psi_\nu\colon \mathbb{R}\to\mathbb{R}$ is locally Lipschitz continuous and is not 
necessarily convex, $f_b\ge 0$ is a constant bound of the magnitude of the friction force.  In particular, 
when $f_b=0$, the last two relations in \eqref{cont6} degenerate to the frictionless condition
\[ -\bsigma_\tau=\bzero \quad{\rm on}\ \Gamma_C. \]
We assume the following properties on the function $\psi_\nu \colon \real \to\real$:
\begin{equation}
\left\{\begin{array}{ll}
{\rm (a)}  \ \psi_\nu(\cdot)\ \mbox{is locally Lipschitz on}\ \real;\\ [0.2mm]
{\rm (b)\  there\ exist\ constants}\ {\bar{c}}_0, {\bar{c}}_1 \ge 0\ {\rm such\ that}\\
{}\qquad  |\partial\psi_\nu(z)| \le {\bar{c}}_0+ {\bar{c}}_1\,|z|\quad \forall\,z\in\real;\\[0.2mm]
{\rm (c)\  there\ exists\ a\ constant}\ \alpha_{\psi_\nu} \ge 0\ {\rm such\ that}\\
{}\qquad \psi_\nu^0(z_1;z_2-z_1)+\psi_\nu^0(z_2; z_1-z_2)\le\alpha_{\psi_\nu}|z_1-z_2|^2\quad \forall\,z_1,z_2\in\real.
\end{array}
\right. \label{psi_nu}
\end{equation}

The weak formulation of the contact problem of \eqref{ve1}--\eqref{ve5} and \eqref{cont6} is the following
VHI (cf.\ \cite{HMS14}, \cite[Problem 4.20]{Han2024}).  

\begin{problem}\label{prob:cont1}
{\it Find a displacement field $\bu\in \bV$ such that}
\begin{align}
& ({\cal E}(\bvarepsilon({\bu})),\bvarepsilon(\bv) -\bvarepsilon(\bu))_\mathbb{Q}
+I_{\Gamma_C}(f_b|\bv_\tau|)-I_{\Gamma_C}(f_b|\bu_\tau|)+I_{\Gamma_C} (\psi_\nu^0 (u_\nu; v_\nu-u_\nu))\nonumber\\
&{}\qquad \ge\langle \fb,\bv-\bu\rangle\quad\forall\,\bv\in \bV.
\label{cont7}
\end{align}
\end{problem}

Well-posedness of Problem \ref{prob:cont1} is stated in Theorem \ref{t21}.  The more general case
of a friction bound $f_b=f_b(u_\nu)$ is considered in \cite{HMS14}, see also \cite[Example 5.51]{Han2024}.

\section{An abstract stationary variational-hemivariational inequality}\label{sec:abs}

Problem \ref{prob:cont1} can be studied within the framework of an abstract stationary VHI discussed in this 
section. Problem \ref{prob:ve} and Problem \ref{prob:VI1} can be studied as special cases of a VI and a VE
of the abstract stationary VHI.  Denote by $\Delta$ the physical domain or its sub-domain, or its boundary 
or part of the boundary, and denote by $I_\Delta$ the integration over $\Delta$,
\[ I_\Delta(v)=\int_\Delta v\,dx\ {\rm if}\ \Delta\subset\Omega,\quad
I_\Delta(v)=\int_\Delta v\,ds\ {\rm if}\ \Delta\subset\Gamma.\]
For a positive integer $m$, we let
\begin{equation}
V_\psi=L^2(\Delta;\mathbb{R}^m).
\label{SpVpsi}
\end{equation}
For application in the study of Problem \ref{prob:cont1}, we take $m=1$; for some other contact problems
(e.g., Problem 4.18 in \cite{Han2024}), we take $m=d$.

We first introduce the following assumptions on the data for the abstract VHI.
\smallskip

\noindent \underline{$H(V)$} $V$ is a real Hilbert space.

\smallskip
\noindent \underline{$H(K)$} $K$ is a non-empty, closed and convex set in $V$.

\smallskip
\noindent \underline{$H(A)$} $A\colon V\to V^*$ is $L_A$-Lipschitz continuous and $m_A$-strongly monotone.

\smallskip
\noindent \underline{$H(\Phi)$} $\Phi\colon V\to\mathbb{R}$ is convex and continuous on $V$.

\smallskip
\noindent \underline{$H(\psi)$}  $\gamma_\psi\in {\cal L}(V;V_\psi)$;  $\psi\colon \mathbb{R}^m\to\mathbb{R}$
is locally Lipschitz continuous and for some non-negative constants $c_\psi$ and $\alpha_\psi$,
\begin{align}
|\partial\psi(z)|_{\mathbb{R}^m} &\le c_\psi\left(1+|z|_{\mathbb{R}^m}\right)
\quad\forall\,z\in\mathbb{R}^m, \label{eq10w} \\
\psi^0(z_1;z_2-z_1)+\psi^0(z_2;z_1-z_2) & \le\alpha_\psi|z_1-z_2|_{\mathbb{R}^m}^2
\quad\forall\,z_1,z_2\in {\mathbb{R}^m}.   \label{eq4}
\end{align}

\noindent \underline{$H(f)$} $f\in V^*$.

\smallskip
Note that an operator $A\colon V\to V^*$ is said to be $L_A$-Lipschitz continuous if
\begin{equation}
\|Av_1-Av_2\|_{V^*} \le L_A \|v_1-v_2\|_V\quad\forall\,v_1,v_2\in V,
\label{eq:5.3b}
\end{equation}
and it is said to be $m_A$-strongly monotone if
\begin{equation}
\langle Av_1-Av_2,v_1-v_2\rangle \ge m_A\|v_1-v_2\|_V^2 \quad\forall\,v_1,v_2\in V.
\label{eq2z}
\end{equation}
A consequence of the assumption $H(\Phi)$ is that $\Phi(\cdot)$ is bounded below by a function that grows
at most linearly, i.e., for some constants $c_3$ and $c_4$, not necessarily positive,
\begin{equation}
\Phi(v)\ge c_3+c_4 \|v\|_V\quad\forall\,v\in V,
\label{Phi:lb}
\end{equation}
cf.\ e.g., \cite[Lemma 11.3.5]{AH2009}.

Generally, we can consider the situation where $\psi=\psi(\bx,z)$ is a function defined for
$\bx\in\Delta$ and $z\in\mathbb{R}^m$.  To simplify the exposition, we will only consider the
case where $\psi=\psi(z)$ does not depend on $\bx\in\Delta$.  We introduce the following assumption.

The condition \eqref{eq4} is equivalent to the following inequality:
\begin{equation}
\langle v^*_1-v^*_2,v_1-v_2\rangle \ge -\alpha_\Psi |v_1-v_2|_{\mathbb{R}^m}^2\quad
\forall\,v_i\in \mathbb{R}^m,\, v^*_i\in\partial \psi(v_i),\ i=1,2.
\label{eq4'}
\end{equation}

It can be shown (e.g., Lemma 4.2 in \cite{HFWH25}) that under the assumption $H(\psi)$,
\begin{equation}
\left|I_\Delta(\psi^0(\gamma_\psi u;\gamma_\psi v))\right|
\le c\left(1+\|u\|_V\right)\|\gamma_\psi v\|_{V_\psi}\quad\forall\,u,v\in V.
\label{3.32a}
\end{equation}

The abstract VHI is the following.

\begin{problem}\label{prob:VHI}
Find $u\in K$ such that
\begin{equation}
\langle Au,v-u\rangle +\Phi(v)-\Phi(u)+I_\Delta(\psi^0(\gamma_\psi u;\gamma_\psi v-\gamma_\psi u))
\ge \langle f,v-u\rangle \quad\forall\,v\in K.
\label{eq1}
\end{equation}
\end{problem}

In \cite{Han2024}, this problem is called a VHI of rank $(1,1)$ to reflect the fact that in the
variational-hemivariational inequality \eqref{eq1}, the convex function $\Phi$ depends on one argument 
and the locally Lipschitz continuous function $\psi$ depends on one argument.  In the general case 
$K\not=V$, Problem \ref{prob:VHI} can be viewed as a constrained VHI of rank $(1,1)$.

When $K=V$ is the entire space, Problem \ref{prob:VHI} becomes an unconstrained VHI of rank $(1,1)$:
Find $u\in V$ such that
\begin{equation}
\langle Au,v-u\rangle +\Phi(v)-\Phi(u)+I_\Delta(\psi^0(\gamma_\psi u;\gamma_\psi v-\gamma_\psi u))
\ge \langle f,v-u\rangle \quad\forall\,v\in V.
\label{eq1aa}
\end{equation}

\subsection{Well-posedness of the abstract stationary variational-hemivariational inequality}\label{subsec:abs}

As an intermediate step in the well-posedness of Problem \ref{prob:VHI}, we analyze an auxiliary VHI.

\begin{problem}\label{prob:aux}
Find $u\in K$ such that
\begin{equation}
\langle Au,v-u\rangle +\Phi(v)-\Phi(u)+\Psi^0(u;v-u)\ge \langle f,v-u\rangle \quad\forall\,v\in K.
\label{eq1aux}
\end{equation}
\end{problem}

In the analysis of Problem \ref{prob:aux}, we assume $H(V)$, $H(K)$, $H(A)$, $H(\Phi)$, $H(f)$, and 

\noindent\underline{$H(\Psi)$} $\Psi\colon V\to\mathbb{R}$ is locally Lipschitz continuous, and for 
some constant $\alpha_\Psi\ge 0$,
\[ \Psi^0(v_1;v_2-v_1)+\Psi^0(v_2;v_1-v_2)\le\alpha_\Psi \|v_1-v_2\|_V^2\quad\forall\,v_1,v_2\in V.\]

\smallskip

In the majority of the references on well-posedness analysis of VHIs, the operator $A\colon V\to V^*$ 
is assumed pseudomonotone, coercive, and strongly monotone, and abstract surjectivity results for 
pseudomonotone operators (e.g., \cite[Theorem 2.12]{NP1995}, \cite[Theorem 1.3.70]{DMP2}) are applied 
(cf.\ \cite{MOS2013, SM2025}).  While such an approach carries its own merit, it is desirable 
to have a more accessible approach for applied and computational mathematicians, and for engineers.
One accessible approach is developed in \cite{Han20, Han21}, and it does not require knowledge on 
pseudomonotone operators or abstract analysis.  We provide a summarizing account of this approach below.
For details, the reader can consult the book \cite{Han2024}.

In the first step of the accessible approach, we assume additionally that $A\colon V\to V^*$ is a
potential operator, i.e., it is the G\^{a}teaux derivative of a functional $F_A\colon V\to\mathbb{R}$.
This allows us to define a related optimization problem.

\begin{problem}\label{prob:opt}
Find $u\in K$ such that
\[ E(u)=\inf\left\{E(v)\mid v\in K\right\} \]
where the energy functional 
\[ E(v)=F_A(v)+\Phi(v)+\Psi(v)-\langle f,v\rangle, \quad v\in V.\]
\end{problem}

Problem \ref{prob:aux} is analyzed through Problem \ref{prob:opt}.  

\begin{theorem}\label{thm-1}
Assume $H(V)$, $H(K)$, $H(A)$, $H(\Phi)$, $H(\Psi)$, $H(f)$, and $\alpha_\Psi<m_A$. Assume additionally 
that $A$ is a potential operator with the potential $F_A$.  Then, Problem \ref{prob:opt} has a unique 
solution $u\in K$, which is also the unique solution of Problem \ref{prob:aux}.
\end{theorem}

In the second step of the accessible approach, we get rid of the additional assumption that $A$ is a
potential operator by a fixed-point technique.  More precisely, for a fixed parameter $\theta>0$, 
given any element $w\in K$, consider the auxiliary problem of finding $u\in K$ such that
\begin{align}
& (u,v-u)_V + \theta \left[ \Phi(v)-\Phi(u)+\Psi^0(u;v-u)\right] \nonumber\\
& {}\qquad \ge (w,v-u)_V - \theta \left[\langle Aw,v-u\rangle- \langle f,v-u\rangle\right]
 \quad\forall\,v\in K. \label{5.54}
\end{align}
An application of Theorem \ref{thm-1} shows that for $\theta>0$ sufficiently small, the inequality
\eqref{5.54} admits a unique solution $u\in K$.  Moreover, it can be shown that the mapping 
$w\mapsto u$ is a contraction.  By the Banach fixed-point theorem, for $\theta>0$ sufficiently small,
the mapping $w\mapsto u$ has a unique fixed-point, and this unique fixed-point is the unique solution of
Problem \ref{prob:aux}.  We state the result next; its detailed proof can be found 
in \cite[Section 5.2]{Han2024}.

\begin{theorem}\label{thm0}
Assume $H(V)$, $H(K)$, $H(A)$, $H(\Phi)$, $H(\Psi)$, $H(f)$, and $\alpha_\Psi<m_A$. 
Then, Problem \ref{prob:aux} has a unique solution $u\in K$.
\end{theorem}

In the last step of the accessible approach, let 
\begin{equation}
\Psi(v)=I_\Delta(\psi(\gamma_\psi v)),\quad v\in V
\label{def:Psi}
\end{equation}
in \eqref{eq1aux}.  Denote by $c_\Delta>0$ the smallest constant in the inequality
\begin{equation}
I_\Delta(|\gamma_\psi v|_{\mathbb{R}^m}^2)\le c_\Delta^2 \|v\|_V^2\quad\forall\,v\in V.
\label{eq12w}
\end{equation}
It can be shown that $H(\psi)$ implies $H(\Psi)$ with $\alpha_\Psi=\alpha_\psi c_\Delta^2$, 
cf.\ Theorem 5.19 in \cite{Han2024}.

In most mathematics references (e.g., \cite{MOS2013, SM2025}), the form of a stationary VHI studied 
is Problem \ref{prob:aux}.  We comment that the inequality \eqref{eq1aux} looks simpler 
than \eqref{eq1}, however, the form \eqref{eq1} arises directly in application problems.  
Any solution of Problem \ref{prob:aux} is a solution of Problem \ref{prob:VHI} since (\cite[Theorem 3.47]{MOS2013})
\[ \Psi^0(u;v)\le I_\Delta(\psi^0(\gamma_\psi u;\gamma_\psi v)). \]
Conversely, if $\psi$ or $-\psi$ is regular in the sense of \eqref{eq:regular}, then 
\[ \Psi^0(u;v)=I_\Delta(\psi^0(\gamma_\psi u;\gamma_\psi v)); \]
as a result, \eqref{eq1aux} and \eqref{eq1} coincide, and Problem \ref{prob:aux} and
Problem \ref{prob:VHI} are identical. This is the approach taken in early references on mathematical
theory of variational-hemivariational inequalities, e.g., \cite{MOS2013}.  

In our approach, we do not assume the regularity of the function $\psi$ or $-\psi$.  
Again, we note that a solution of Problem \ref{prob:aux} is also
a solution of Problem \ref{prob:VHI}.  In addition, it can be shown that a solution of 
Problem \ref{prob:VHI} is unique.  Thus, Problem \ref{prob:VHI} admits a unique solution.
The result is stated next, cf.\ \cite[Section 5.4]{Han2024} for details.

\begin{theorem}\label{thm1}
Assume $H(V)$, $H(K)$, $H(A)$, $H(\Phi)$, $H(\psi)$, $H(f)$, and $\alpha_\psi c_\Delta^2<m_A$.  
Then Problem \ref{prob:VHI} has a unique solution $u\in K$. Moreover, the solution $u\in K$ 
depends Lipschitz continuously on $f\in V^*$.
\end{theorem}

\subsection{Galerkin method for the abstract stationary variational-hemivariational inequality}

Since there is no analytic solution formula for a VHI arising in applications, numerical methods are 
needed to solve the problem.  In this subsection, we discuss the numerical solution of Problem \ref{prob:VHI}.
The numerical method is of Galerkin type.  We present a convergence result for the numerical solutions  
and a C\'{e}a-type inequality for error estimation of the numerical solutions.  For Problem \ref{prob:VHI}, 
we make the assumptions stated in Theorem \ref{thm1}, so as to guarantee that the problem has a unique solution.
Let $V^h$ be a finite dimensional subspace of $V$, $h>0$ being a spatial discretization parameter.
Let $K^h$ be a non-empty, closed and convex subset of $V^h$.
Then, a Galerkin approximation of Problem~\ref{prob:VHI} is the following.

\smallskip
\begin{problem} \label{prob:VHIh}
{\it Find an element $u^h\in K^h$ such that}
\begin{align}
& \langle Au^h,v^h-u^h\rangle+\Phi(v^h)-\Phi(u^h)+I_\Delta(\psi^0(\gamma_\psi u^h;\gamma_\psi v^h-\gamma_\psi u^h))\nonumber \\
&\qquad{}\ge\langle f,v^h-u^h\rangle\quad\forall\,v^h\in K^h.
\label{eq1h}
\end{align}
\end{problem}

For the well-posedness of Problem \ref{prob:VHIh}, we can apply Theorem \ref{thm1} which is valid in the  
setting of finite-dimensional spaces as well.  For completeness, we state the result formally as a theorem.

\begin{theorem}\label{thm1h}
Keep the assumptions stated in Theorem \ref{thm1}.  Let $V^h$ be a finite-dimensional subspace of $V$ and 
let $K^h$ be a non-empty, closed and convex subset of $V^h$. Then Problem~\ref{prob:VHIh} has a unique solution.
\end{theorem}

The approximation is called external if $K^h\not\subset K$, and is internal if $K^h\subset K$.  In 
\cite{HSD18}, the internal approximation with the choice $K^h=V^h\cap K$ is considered for Problem \ref{prob:VHI}.

The following convergence result is proved in \cite[Section 4.3]{HFWH25} and it follows \cite{HZ19}.  
An important point about this convergence result is that we do not assume any solution regularity 
other than the basic regularity $u\in V$ guaranteed in the well-posedness result, namely, Theorem \ref{thm1}.

\smallskip
\begin{theorem}\label{thm:converge}
Keep the assumptions stated in Theorem \ref{thm1}.  Moreover, assume $V^h$ is a finite-dimensional subspace of $V$,
$K^h$ is a non-empty, closed and convex subset of $V^h$, and 
\begin{align}
&\label{3.3i} v^h\in K^h\ {\rm and}\ v^h\weak v\ {\rm in}\ V\ {\rm imply}\ v\in K;\\
&\label{3.3h} \forall\,v\in K,\ \exists\,v^h\in K^h\ {\rm such\ that\ }v^h\to v\ {\rm in\ }V\ {\rm as}\ h\to 0.
\end{align}
Let $u$ and $u^h$ be the solutions of Problem \ref{prob:VHI} and Problem \ref{prob:VHIh}, respectively.  Then,
\begin{equation}
 u^h\to u\quad{\rm in\ }V\ {\rm as\ }h\to 0. \label{4.9a}
\end{equation}
\end{theorem}

Theorem \ref{thm1} and Theorem \ref{thm:converge} are rather general, and here we consider two special cases.

First, we consider the case of a HVI with the choice $\Phi\equiv 0$ in Problem \ref{prob:VHI}.

\begin{problem}\label{P21}
{\it Find an element  $u \in K$ such that}
\begin{equation}
\langle Au,v-u\rangle+I_\Delta(\psi^0(\gamma_\psi u;\gamma_\psi v-\gamma_\psi u))\ge\langle f,v-u\rangle\quad\forall\,v\in K.
\label{hv2}
\end{equation}
\end{problem}

The corresponding numerical method Problem \ref{prob:VHIh} takes the following form.

\begin{problem} \label{P2h}
{\it Find an element $u^h\in K^h$ such that}
\begin{equation}
\langle Au^h,v^h-u^h\rangle+I_\Delta(\psi^0(\gamma_\psi u^h;\gamma_\psi v^h-\gamma_\psi u^h))
\ge\langle f,v^h-u^h\rangle\quad\forall\,v^h\in K^h.
\label{hvh2}
\end{equation}
\end{problem}

\begin{theorem}\label{thm:converge2}
Assume $H(V)$, $H(K)$, $H(A)$, $H(\psi)$,  $H(f)$, and $\alpha_\psi c_\Delta^2<m_A$.
Moreover, assume $V^h$ is a finite-dimensional subspace of $V$, $K^h$ is a non-empty, closed and convex subset
of $V^h$, and \eqref{3.3i}--\eqref{3.3h} hold.  Then, Problem \ref{P21} admits a unique solution $u\in K$, 
Problem \ref{P2h} admits a unique solution $u^h\in K^h$, and we have the convergence:
\[ u^h\to u\quad{\rm in\ }V\ {\rm as\ }h\to 0. \]
\end{theorem}

As another particular case, we consider a VI, obtained from Problem \ref{prob:VHI}
by setting $\psi\equiv 0$.

\begin{problem}\label{P4}
{\it Find an element  $u \in K$ such that}
\begin{equation}
\langle Au,v-u\rangle+\Phi(v)-\Phi(u)\ge\langle f,v-u\rangle\quad\forall\,v\in K.
\label{hv4}
\end{equation}
\end{problem}

The numerical method is the following.

\begin{problem} \label{P4h}
{\it Find an element $u^h\in K^h$ such that}
\begin{equation}
\langle Au^h, v^h - u^h \rangle + \Phi(v^h)- \Phi(u^h)\ge\langle f,v^h-u^h\rangle\quad\forall\,v^h\in K^h.
\label{hvh4}
\end{equation}
\end{problem}

\begin{theorem}\label{thm:converge4}
Assume $H(V)$, $H(K)$, $H(A)$, $H(\Phi)$, and $H(f)$.  Moreover, assume $V^h$ is a finite-dimensional subspace
of $V$, $K^h$ is a non-empty, closed and convex subset of $V^h$, and \eqref{3.3i}--\eqref{3.3h} hold.
Then, Problem \ref{P4} admits a unique solution $u\in K$, Problem \ref{P4h} admits a unique solution 
$u^h\in K^h$, and we have the convergence:
\[ u^h\to u\quad{\rm in\ }V\ {\rm as\ }h\to 0. \]
\end{theorem}

For error estimates of the numerical solution defined by Problem \ref{prob:VHIh} for the approximation of 
the solution of Problem \ref{prob:VHI}, let $v\in K$ and $v^h\in K^h$ be arbitrary and define
\begin{equation}
R_u(v,w):=\langle Au,v-w\rangle+\Phi(v) -\Phi(w) + I_\Delta(\psi^0(\gamma_\psi u;\gamma_\psi v-\gamma_\psi w))
-\langle f,v-w\rangle, \label{err8} 
\end{equation}
The following result is proved as Theorem 7 in \cite{HFWH25}.

\begin{theorem}\label{thm:num_1}
Assume $H(K)$, $H(A)$, $H(\Phi)$, $H(\psi)$,  $H(f)$, and $\alpha_\psi c_\Delta^2<m_A$.
Then for the solution $u$ of Problem \ref{prob:VHI} and the solution $u^h$ of Problem \ref{prob:VHIh},
we have the C\'{e}a-type inequality
\begin{equation}
\|u-u^h\|_V^2\le c\inf_{v^h\in K^h}\left[\|u-v^h\|_V^2+\|\gamma_\psi (u-v^h)\|_{V_\psi}+R_u(v^h,u)\right]
+c\inf_{v\in K} R_u(v,u^h). \label{err16}
\end{equation}
\end{theorem}

In the case of an internal approximation, $K^h\subset K$, and
\[ \inf_{v\in K} R_u(v,u^h) =0. \]
Then, the C\'{e}a-type inequality \eqref{err16} simplifies to
\begin{equation}
\|u-u^h\|_V^2 \le  c\inf_{v^h\in K^h}\left[\|u-v^h\|_V^2 +\|\gamma_\psi (u-v^h)\|_{V_\psi}+R_u(v^h,u)\right].
\label{err16.aa}
\end{equation}

We note that \eqref{err16} and \eqref{err16.aa} are generalizations of the C\'{e}a-type inequality from 
the numerical solution of VIs (cf.\ \cite{Fa74, KO1988, HS2002}) to that of VHIs.  To proceed further, 
we need to bound the residual term \eqref{err8} and this depends on the problem to be solved.

\section{Studies of the sample variational-hemivariational inequality}\label{sec:contact}

In this section, we return to Problem~\ref{prob:cont1} by applying the theoretical results reviewed in 
Section \ref{sec:abs}.  We first explore the solution existence and uniqueness, then introduce a 
linear finite element method to solve the problem and present an optimal order error estimate 
under certain solution regularity assumptions.  Details can be found in \cite{HFWH25}.

Let $\lambda_\nu>0$ be the smallest eigenvalue of the eigenvalue problem
\[ \bu\in \bV,\quad \int_\Omega \bvarepsilon({\bu}):\bvarepsilon(\bv)\,dx
=\lambda \int_{\Gamma_C} u_\nu v_\nu \, ds\quad\forall\,\bv\in \bV.\]
Then, we have the trace inequality
\[ \|v_\nu\|_{L^2(\Gamma_C)} \le \lambda_\nu^{-1/2} \|\bv\|_{\boldsymbol V}\quad\forall\,\bv\in \bV.\]
By applying Theorem~\ref{thm1}, an existence and uniqueness result can be established for Problem~\ref{prob:cont1}.

\begin{theorem}\label{t21}
Assume \eqref{Ass:E}, \eqref{Ass:f}, \eqref{psi_nu}, $f_b\ge 0$, and
$\alpha_{\psi_\nu} \lambda_{\nu}^{-1} < m_{{\cal E}}$.  Then Problem \ref{prob:cont1} has a unique solution $\bu\in \bV$.
\end{theorem}

Theorem~\ref{t21} provides the existence of a unique displacement field $\bu\in \bV$ of the contact problem.
The stress field $\bsigma\in \mathbb{Q}$ is uniquely determined by using the constitutive law \eqref{ve2}.

Then, we turn to the discretization of Problem \ref{prob:cont1} using the finite element method.
For simplicity, assume $\Omega$ is a polygonal domain ($d=2$) or a polyhedral domain ($d=3$) .  
We express the three parts $\Gamma_D$, $\Gamma_N$, $\Gamma_C$ of the boundary
as unions of closed flat components with disjoint interiors:
\[ \overline{\Gamma_Z}=\cup_{i=1}^{i_Z}\Gamma_{Z,i},\quad Z=D,N,C.\]
Let $\{{\cal T}^h\}$ be a regular family of partitions of $\overline{\Omega}$
into triangles ($d=2$) or tetrahedrons ($d=3$) that are compatible with the partition of
the boundary $\partial\Omega$ into $\Gamma_{Z,i}$, $1\le i\le i_Z$, $Z=D,N,C$,
in the sense that if the intersection of one side/face of an
element with one set $\Gamma_{Z,i}$ has a positive measure with
respect to $\Gamma_{Z,i}$, then the side ($d=2$) or the face ($d=3$) lies entirely in $\Gamma_{Z,i}$.
The corresponding linear finite element space is
\begin{equation}
\bV^h=\left\{\bv^h\in C(\overline{\Omega})^d \mid \bv^h|_T\in \mathbb{P}_1(T)^d
   \ {\rm for}\ T\in {\cal T}^h,\ \bv^h=\bzero\ {\rm on\ }\overline{\Gamma_D}\right\}.
\label{Vh}
\end{equation}

The finite element approximation of Problem \ref{prob:cont1} is the following.

\begin{problem}\label{p2h}
{\it Find a displacement field $\bu^h\in \bV^h$ such that}
\begin{align}
&\left({\cal E}(\bvarepsilon(\bu^h)),\bvarepsilon(\bv^h)-\bvarepsilon(\bu^h)\right)_{\mathbb{Q}}
+ I_{\Gamma_C}(f_b\,|\bv^h_\tau|) - I_{\Gamma_C}(f_b\,|\bu^h_\tau|) \nonumber\\[1mm]
&\qquad{}+I_{\Gamma_C}(\psi^0_\nu(u_\nu^h;v_\nu^h-u_\nu^h))\ge \langle\fb,\bv^h-\bu^h\rangle\quad\forall\,\bv^h\in \bV^h.
\label{s7.2a}
\end{align}
\end{problem}

Similar to Theorem \ref{t21}, we know that Problem \ref{p2h} admits a unique solution $\bu^h\in \bV^h$.
Convergence of the finite element solution follows from Theorem \ref{thm:converge}:
\[ \|\bu^h-\bu\|_{\boldsymbol V}\to 0\quad {\rm as}\ h\to 0. \]
We comment that the two conditions \eqref{3.3i}--\eqref{3.3h} are valid with $K=\bV$, $K^h=\bV^h$ for the
finite element space \eqref{Vh}.

For an error analysis, we start with
\begin{equation}
\|\bu-\bu^h\|_{\boldsymbol V}^2  \le c\inf_{{\boldsymbol v}^h\in {\boldsymbol V}^h}\left[\|\bu-\bv^h\|_{\boldsymbol V}^2
+\|u_\nu-v^h_\nu\|_{L^2(\Gamma_C)}+R_{\boldsymbol u}(\bv^h,\bu)\right],
\label{s7.2aaa}
\end{equation}
where the residual-type term as defined in \eqref{err8} is
\begin{align}
R_{\boldsymbol u}(\bv^h,\bu)&=({\cal E}(\bvarepsilon(\bu)),\bvarepsilon(\bv^h-\bu))_{\mathbb{Q}}
+  I_{\Gamma_C}(f_b\left(|\bv_\tau^h|-|\bu_\tau|\right)) \nonumber\\[1mm]
&\quad{}+I_{\Gamma_C} (\psi^0_\nu (u_\nu; v_\nu^h-u_\nu))- \langle\fb, \bv^h-\bu\rangle. \label{s7.2c}
\end{align}

To derive an error bound, we need to make solution regularity assumptions:
\begin{equation}
\bu\in H^2(\Omega;\mathbb{R}^d),\quad \bsigma={\cal E}(\bvarepsilon(\bu))\in H^1(\Omega;\mathbb{S}^d),
\quad  \bu|_{\Gamma_{C,i}}\in H^2(\Gamma_{C,i};\real^d), \quad 1\le i\le i_C.
\label{s7.2b}
\end{equation}
For homogeneous materials, ${\cal E}$ does not depend on the spatial variable $\bx$ and the second regularity 
is a consequence of the first one in \eqref{s7.2b}.  Under the first two regularity assumptions in \eqref{s7.2b},
it can be shown that the pointwise relations hold:
\begin{align}
&{\rm div}\,{\cal E}(\bvarepsilon({\bu}))+\fb_0=\bzero\quad{\rm a.e.\ in\ }\Omega, \label{s7.2d}\\
&\bsigma\bnu=\fb_2\quad{\rm a.e.\ on\ }\Gamma_N. \label{s7.2e}
\end{align}
Then,
\[ R_{\boldsymbol u}(\bv^h,\bu)=\int_{\Gamma_C} \left[ \bsigma\bnu{\cdot}(\bv^h-\bu)
 +f_b \left(|\bv^h_\tau|-|\bu_\tau|\right)+\psi_\nu^0(u_\nu;v^h_\nu-u_\nu)\right] ds,\]
and we can bound
\begin{equation}
\left|R_{\boldsymbol u}(\bv^h,\bu)\right| \le c\,\|\bu-\bv^h\|_{L^2(\Gamma_C;\mathbb{R}^d)}.
\label{s7.2g}
\end{equation}
So, from \eqref{s7.2aaa}, we can derive the C\'{e}a-type inequality
\begin{equation}
\|\bu-\bu^h\|_{\boldsymbol V} \le c\inf_{{\boldsymbol v}^h\in {\boldsymbol V}^h}\left[\|\bu-\bv^h\|_{\boldsymbol V}
+\|\bu-\bv^h\|_{L^2(\Gamma_C;\mathbb{R}^d)}^{1/2}\right].
\label{s7.2h}
\end{equation}

Due to the first regularity condition in \eqref{s7.2b}, by the Sobolev embedding 
$H^2(\Omega)\subset C(\overline{\Omega})$ valid for $d\le 3$, we know that 
$\bu\in C(\overline{\Omega};\mathbb{R}^d)$
and so its finite element interpolant $\Pi^h\bu\in \bV^h$ is well defined.  Moreover, the following 
error estimate holds (cf.\ \cite{AH2009, BS2008, Ci1978}): for some constant $c>0$ independent of $h$,
\begin{equation}
\|\bu-\Pi^h\bu\|_{L^2(\Omega;\mathbb{R}^d)}+h\,\|\bu-\Pi^h\bu\|_{H^1(\Omega;\mathbb{R}^d)}
\le c\,|\bu|_{H^2(\Omega;\mathbb{R}^d)},
\label{estimate}
\end{equation}
and
\begin{equation}
\|\bu-\Pi^h\bu\|_{L^2(\Gamma_C;\mathbb{R}^d)}\le c\,h^2.
\label{5.20a}
\end{equation}
Then we derive from \eqref{s7.2h} the following optimal order error bound
\begin{equation}
\|\bu-\bu^h\|_{\boldsymbol V}\le c\left[\|\bu-\Pi^h\bu\|_{\boldsymbol V}
+\|\bu-\Pi^h\bu\|_{L^2(\Gamma_C;\mathbb{R}^d)}^{1/2}\right]\le c\,h,
\label{s7.2j}
\end{equation}
where the constant $c$ depends on the quantities $\|\bu\|_{H^2(\Omega;\mathbb{R}^d)}$,
$\|\bsigma\bnu\|_{L^2(\Gamma_C;\mathbb{R}^d)}$ and $\|\bu\|_{H^2(\Gamma_{C,i};\mathbb{R}^d)}$
for $1\le i\le i_C$.

The reader is referred to \cite{HFWH25} for numerical examples providing numerical convergence orders
that match the error estimate \eqref{s7.2j}.

\smallskip
We comment that similar results hold for Problem \ref{prob:VI1} and Problem \ref{prob:ve} by applying
Theorem \ref{thm1} and Theorem \ref{thm:num_1} for the special cases of a VI and a VE.

\section{Mixed variational-hemivariational inequalities in fluid mechanics}\label{sec:fluid}

In the previous sections, we considered VHIs from contact mechanics.  In this section, we consider 
sample VHIs in fluid mechanics.  Since Fujita's pioneering work \cite{Fu93, Fu94}, there has been 
steady progress on the modeling, mathematical analysis and numerical approximation of boundary or 
initial-boundary value problems for flows of viscous incompressible fluid involving nonsmooth slip
or leak boundary conditions. When the nonsmooth boundary condition is of monotone type, 
the mathematical formulation of the problem is a mixed variational inequality.  A large number of 
papers have been published on well-posedness analysis and the numerical solution of Stokes and 
Navier-Stokes variational inequalities.  The reader may consult \cite[Chapter 8]{Han2024} for some 
sample references on this topic. When the nonsmooth boundary condition is allowed to be of
nonmonotone type, the weak formulation of the problem is a mixed hemivariational inequality. 
Various VHIs arising in fluid mechanics have been studied in the literature, e.g.,
\cite{CMMY23, HN23, MCHD24, MD22, TC25}.  In this section, we present two mathematical models:
a Stokes HVI and a Navier-Stokes HVI.

Let $\Omega\subset \mathbb{R}^d$ ($d\leq 3$ in applications) be a Lipschitz domain.
Its boundary $\Gamma=\partial\Omega$ is decomposed into two parts:
$\Gamma=\overline{\Gamma_D}\cup\overline{\Gamma_S}$ such that $|\Gamma_D|>0$, 
$|\Gamma_S|>0$, and $\Gamma_D\cap\Gamma_S=\emptyset$. We will impose a Dirichlet boundary 
condition on $\Gamma_D$ and a slip boundary condition of friction type on $\Gamma_S$.  Denote by 
$\bnu$ the unit outward normal to $\Gamma$. For a vector-valued function $\bu$ on the boundary, 
let $u_\nu=\bu\cdot\bnu$ and $\bu_{\tau}=\bu-u_\nu\bnu$ be the normal component and the tangential
component, respectively. Let $\mu>0$ be the viscosity coefficient, and let $\fb$ be the source function.

The pointwise relations of a sample Stokes HVI are
\begin{align}
& -{\rm div}(2\,\mu\,\bvarepsilon(\bu))+\nabla p=\fb\quad{\rm in}\ \Omega,
\label{eq:8.56e}\\
&\, \textrm{div}\,\bu=0\quad{\rm in}\ \Omega,
\label{eq:8.57e}\\
& \bu=\bzero\quad{\rm on}\ \Gamma_D, \label{eq:8.58e}\\
& u_\nu=0,\quad-\bsigma_\tau\in \partial \psi_\tau(\bu_\tau)\ {\rm on}\ \Gamma_S.  \label{eq:8.59e}
\end{align}
The condition $u_\nu=0$ in \eqref{eq:8.59e} means that the fluid can not pass through $\Gamma_S$ 
outside the domain, and this condition is usually called the no-leak condition. The second part 
in \eqref{eq:8.59e}  represents a friction condition, relating the frictional force $\bsigma_{\tau}$ 
with the tangential velocity $\bu_\tau$.  In \eqref{eq:8.59e}, the function 
$\psi_\tau\colon \mathbb{R}^d\to\mathbb{R}$ is locally
Lipschitz continuous.  We assume the following properties on the function $\psi_\tau$:
\begin{equation}
\left\{\begin{array}{ll}
{\rm (a)}  \ \psi_\tau(\cdot)\ \mbox{is locally Lipschitz on}\ \real^d;\\ [0.2mm]
{\rm (b)\  there\ exist\ constants}\ {\bar{c}}_0, {\bar{c}}_1 \ge 0\ {\rm such\ that}\\
{}\qquad  |\partial\psi_\tau(\bz)| \le {\bar{c}}_0+ {\bar{c}}_1\,|\bz|\quad \forall\,\bz\in\real^d;\\[0.2mm]
{\rm (c)\  there\ exists\ a\ constant}\ \alpha_{\psi_\tau} \ge 0\ {\rm such\ that}\\
{}\qquad \psi_\tau^0(\bz_1;\bz_2-\bz_1)+\psi_\tau^0(\bz_2; \bz_1-\bz_2)\le\alpha_{\psi_\tau}|\bz_1-\bz_2|^2\quad \forall\,\bz_1,\bz_2\in\real^d.
\end{array}
\right. \label{6.4a}
\end{equation}

We comment that it is possible to consider the more general case $\psi_\tau\colon \Gamma_S\times
\mathbb{R}^d\to\mathbb{R}$ where $\psi_\tau$ depends on the spatial variable $\bx$.

\smallskip
In this section, we define the function space
\begin{equation}
\bV=\left\{\bv\in\bH^1(\Omega)\mid\bv=\bzero\ {\rm on}\ \Gamma_D,\ v_\nu=0\ {\rm on}\ \Gamma_S\right\},
\label{Sp:Vf}
\end{equation}
for the velocity field, where $\bH^1(\Omega)=H^1(\Omega;\mathbb{R}^d)$. As a consequence of Korn's inequality, the quantity
\[ \|\bvarepsilon(\bv)\|_{L^2(\Omega;\mathbb{S}^d)}:=
\left(\int_\Omega|\bvarepsilon(\bv)|^2\,dx\right)^\frac{1}{2}\]
defines a norm for $\bv\in \bV$ and it is equivalent to the standard $\bH^1(\Omega)$-norm on $\bV$. We use
\begin{equation}
\|\cdot\|_{\boldsymbol V}=\|\bvarepsilon(\cdot)\|_{L^2(\Omega;\mathbb{S}^d)}
\label{eq:8.12e}
\end{equation}
for the norm on $\bV$.  For the pressure variable, we use the space
\begin{equation}
Q=L^2_0(\Omega)=\left\{q\in L^2(\Omega)\mid I_\Omega(q)= 0\right\}.
\label{Sp:Qf}
\end{equation}
This is a Hilbert space with the standard inner product
\[ (p,q)_{0,\Omega}:=\int_\Omega p\,q\,dx \]
and the corresponding $\|\cdot\|_{0,\Omega}$-norm.

Define the following forms:
\begin{align}
& a(\bu,\bv)=2\mu\int_\Omega \bvarepsilon(\bu):\bvarepsilon(\bv)\,dx\quad\forall\,\bu,\bv\in\bV,\label{f:a}\\
& b(\bv,q)=\int_\Omega q\,{\rm div}\bv\, dx\quad\forall\,\bv\in\bV,\,q\in Q. \label{f:b}
\end{align}
The weak formulation of the problem \eqref{eq:8.56e}--\eqref{eq:8.59e} is derived with the standard approach.

\begin{problem}\label{p:hv01}\index{Stokes HVI}
Find $(\bu,p)\in \bV\times Q$ such that
\begin{align}
&a(\bu,\bv)-b(\bv,p)+I_{\Gamma_S}(\psi_\tau^0(\bu_\tau;\bv_\tau))\geq\langle\fb,\bv\rangle
\quad\forall\,\bv\in \bV,\label{dp05}\\
&b(\bu,q)=0\quad\forall\,q\in Q.\label{dp06}
\end{align}
\end{problem}

Let $\lambda_\tau>0$ be the smallest eigenvalue of the eigenvalue problem
\begin{equation}
\bu\in \bV,\quad \int_\Omega \bvarepsilon({\bu}):\bvarepsilon(\bv)\,dx
=\lambda \int_{\Gamma_S} \bu_\tau \cdot \bv_\tau \, ds\quad\forall\,\bv\in \bV. 
\label{6.13a}
\end{equation}
Then, we have the trace inequality
\[ \|\bv_\tau\|_{L^2(\Gamma_S;\mathbb{R}^d)}\le\lambda_\tau^{-1/2} \|\bv\|_{\boldsymbol V}\quad\forall\,\bv\in \bV.\]
The following well-posedness result is labelled as Theorem 8.22 in \cite{Han2024}.

\begin{theorem}\label{thm:S_HVI}
Assume \eqref{6.4a} and $\alpha_{\psi_\tau}< 2\mu\lambda_\tau$.  Then, for any $\fb\in \bV^*$, Problem \ref{p:hv01} 
has a unique solution $(\bu,p)\in \bV\times Q$ which depends Lipschitz continuously on $\fb\in \bV^*$.
\end{theorem}

We skip details of discussions of the numerical solution of Problem \ref{p:hv01} in this paper, and refer the
reader to \cite{FCHCD20} on a mixed finite element method and to \cite{LHZ22} on stabilized mixed finite 
element methods to solve Problem \ref{p:hv01}.

\smallskip

We now turn to a sample Navier-Stokes HVI.  The pointwise formulation of the problem is
\begin{align}
& -{\rm div}(2\mu\,\bvarepsilon(\bu))+(\bu{\cdot}\nabla)\bu+\nabla p=\fb\quad\textrm{in}\ \Omega, \label{ns01}\\
&\, \textrm{div}\,\bu=0\quad\textrm{in}\ \Omega,  \label{ns02}\\
&\bu=\pmb 0\quad\textrm{on}\ \Gamma_D, \label{ns03}\\
&u_\nu=0,\quad-\bsigma_\tau\in \partial \psi_\tau(\bu_\tau)\quad\textrm{on}\ \Gamma_S.  \label{ns04}
\end{align}
Use the function space $\bV$ defined in \eqref{Sp:Vf} for the velocity, and the function space $Q$ 
defined in \eqref{Sp:Qf} for the pressure.  Use the forms $a(\bu,\bv)$ defined in \eqref{f:a}, 
$b(\bv,q)$ defined in \eqref{f:b}.  Furthermore, define 
\begin{equation}
d(\bu,\bv,\bw) = \int_\Omega (\bu{\cdot}\nabla)\bv {\cdot} \bw\,dx\quad\forall\,\bu,\bv, \bw\in\bV. \label{f:d}
\end{equation}
We continue to assume $\fb\in\bV^*$ and \eqref{6.4a}.  Then, the weak formulation of the problem is

\begin{problem}\label{NS-HVI}\index{Navier-Stokes HVI}
Find $\bu\in \bV$ and $p\in Q$ such that
\begin{align}
&a(\bu,\bv)+d(\bu,\bu,\bv)-b(\bv,p)+I_{\Gamma_S}(\psi_\tau^0(\bu_\tau;\bv_\tau))
\ge\langle\fb,\bv\rangle\quad\forall\,\bv\in \bV,
\label{ns15}\\
&b(\bu,q) = 0 \quad\forall\,q\in Q.
\label{ns16}
\end{align}
\end{problem}

The following well-posedness result on Problem \ref{NS-HVI} is shown in \cite[Section 8.33]{Han2024}.

\begin{theorem}\label{thm:ns}
Assume \eqref{6.4a} and 
\[ \alpha_{\psi_\tau} \lambda_\tau^{-1}<2\,\mu,\quad 
\alpha_{\psi_\tau} \lambda_\tau^{-1}+ c_d M_{\boldsymbol f}<2\,\mu, \]
where $c_d>0$ is a constant in the boundedness inequality
\[ \left|d(\bu,\bv,\bw)\right|\le c_d\|\bu\|_{\boldsymbol V}\|\bv\|_{\boldsymbol V}\|\bw\|_{\boldsymbol V}
\quad\forall\,\bu,\bv,\bw\in \bV,\]
the constant $M_{\boldsymbol f}$ is defined by
\[ M_{\boldsymbol f}=\frac{c_0 \lambda_\tau^{-1/2} |\Gamma_S|^{1/2}
+\|\fb\|_{\boldsymbol V^*}}{2\,\mu-\alpha_{\psi_\tau} \lambda_\tau^{-1}} \]
and $\lambda_\tau>0$ is the smallest eigenvalue of the eigenvalue problem \eqref{6.13a}. Then, 
Problem \ref{NS-HVI} admits a unique solution $(\bu,p)\in \bV\times Q$, 
$\|\bu\|_{\boldsymbol V}\le M_{\boldsymbol f}$, and 
$(\bu,p)\in \bV\times Q$ depends Lipschitz continuously on $\fb\in \bV^*$.
\end{theorem}

Once again, we skip a detailed discussion of the numerical solution of Problem \ref{NS-HVI} in this paper,
and refer the reader to \cite{HCJ21} on a mixed finite element method and to \cite{HJY23} on 
stabilized mixed finite element methods to solve Problem \ref{NS-HVI}.

Finally, we mention that research on VHIs for problems in fluid mechanics is not limited to the 
standard Stokes and Navier-Stokes quations.  E.g., HVIs are analyzed together with their finite element 
solutions for incompressible fluid flows with damping in \cite{HQM24} for the Stokes equations,
and in \cite{WCH26} for the Navier–Stokes equations.  In \cite{AM25}, a HVI is studied for 
incompressible fluid flows with damping and pumping effects.

\section{Miscellaneous remarks}\label{sec:other}

In the previous sections, we discussed the finite element method to solve stationary and mixed VHIs. 
Other numerical methods have been applied to solve VHIs as well.  E.g., the discontinuous Galerkin 
method is applied to solve a HVI for semipermeable media in \cite{WQ20} and to solve a contact problem 
in \cite{WSW23}.  Take the virtual element method (VEM) as another example.  The VEM was first proposed 
and analyzed in \cite{beirao2013basic, beirao2013virtual}. The method has since been applied to a wide 
variety of mathematical models from applications in science and engineering thanks to its strengths 
in handling complex geometries and problems requiring high-regularity solutions. The VEM was first 
applied to solve contact problems in \cite{wriggers2016virtual}.  Further applications of the VEM to 
solve contact problems can be found in a number of publications, e.g., \cite{AHABW20, CHKW22, WZ21, WWH24}.
The VEM was first applied to solve a HVI in \cite{FHH19}. Further references on application of the VEM 
to solve VHIs are \cite{FHH21a, FHH21b, FHH22, LWH20, WWH21, WWH22, XL23a, XL23b}.

In \cite{BHM19}, a numerical method is applied to solve a VHI modeling contact problems for locking 
materials.  An error estimate is derived in the paper; however, derivation of an optimal order error 
bound remains open.

To solve discretized VHIs, optimization based numerical algorithms are developed in \cite{JO21} and
\cite{OJB22} for stationary VHIs, and in \cite{JOB23} for time-dependent VHIs, all with applications to 
contact mechanics.  A theoretical foundation for the optimization approach is explored in \cite{Han20},
where it is shown that certain special VHIs are equivalent to optimization problems.  The possibility of
reformulating some VHIs as optimization problems is a starting point to develop deep neural networks 
to solve those VHIs (\cite{HWW22}).  Deep neural network techniques have also been explored to solve 
an obstacle problem in \cite{CSWL23}.  For another kind of neural network methods to solve obstacle 
problems, see \cite{WD24}.  More research is expected on neural network methods to solve VHIs.

Most of the papers on VHIs deal with PDEs of second-order.  In \cite{FHH22}, a nonconforming virtual element method is studied for a fourth-order HVI in the Kirchhoff plate problem. In \cite{QWLZ23}, 
an interior penalty virtual element method is developed to solve a fourth-order HVI.  In \cite{QLWW25}, 
a nonconforming finite element method is developed to solve a fourth-order history-dependent HVI.  

Although most of the papers on the numerical solution of VHIs deal with the stationary 
(i.e., time-independent) case, various papers are available on numerical methods for solving 
time-dependent VHIs.  Some representative publications are as follows.  In \cite{BBHJ15}, a hyperbolic HVI 
arising in a dynamic contact problem is considered.  Both spatially semi-discrete schemes and fully
discrete schemes are introduced and analyzed.  Optimal order error estimates are derived with the use
of linear finite elements. In \cite{BBH17}, numerical analysis is provided on a parabolic VHI with an 
application to a contact problem for viscoelastic bodies.  In \cite{HH19}, numerical analysis is 
provided on another parabolic VHI with an application to a dynamic contact problem for viscoelastic 
bodies. In \cite{Ba20}, a fully discrete scheme is applied to solve a coupled
system of HVIs for a nonmonotone dynamic contact problem of a non-clamped piezoelectric viscoelastic body.

The term history-dependent (quasi-)variational inequalities first appears in the paper \cite{SM11}.
History-dependent VHIs are first analyzed in \cite{SM16}.  The numerical solution of history-dependent 
VHIs by the finite element method for spatial discretizations is studied in \cite{XHHCW19, XHHCW20, WXHC20}. 
The virtual element method is applied to solve a history-dependent VHI in \cite{XL23}, and is applied to solve a history-dependent mixed VHI in \cite{LXH24}, both concerning applications in contact problems.  

In \cite{XCY25}, numerical analysis for a system of fractional differential hemivariational inequalities 
in the consideration of a thermoviscoelastic frictional contact problem involving time-fractional order 
operators and long memory effects.  In \cite{SCX24}, numerical analysis for a coupled system of a
dynamic HVI, a parabolic VI and a parabolic PDE that describes a thermoviscoelastic contact problem with 
damage and long memory.


Besides contact mechanics and fluid mechanics, HVIs appear also in other fields of applications.  
As an example, a HVI from Bean's model for an irreversible and hysteretic magnetization process of 
high-temperature superconductors in a magnetic field is studied in \cite{HLW23}.

This paper provides an introduction of basic mathematical theory and numerical analysis of VHIs, 
focusing on the stationary ones. Well-posedness analysis of VHIs is conducted with an approach 
that does not require knowledge of abstract theory of pseudomonotone operators.  The main idea of
this accessible approach for the study of stationary VHIs is explained in the paper.  
We hope this helps to attract more researchers in applied and computational mathematics, engineering 
to the exciting and challenging area of VHIs.  

A basic theory is sketched for the numerical solution of stationary VHIs.  The numerical scheme is 
constructed within the framework of Galerkin methods, and in particular, by the finite element method.  
Convergence of the numerical solutions is explored, and derivation of error estimates for the numerical 
solutions is demonstrated through the approximation of a stationary VHI in contact mechanics.  
Besides the finite element method, references are also provided for the use of the virtual element method 
and the discontinuous Galerkin method. Further efforts are needed for the development of efficient and 
effective solution algorithms to solve the discretized VHIs and for the analysis of such algorithms.  
It is even possible to develop machine learning methods to solve VHIs.


\begin{thebibliography}{99}

\bibitem{AF2003}
R. A. Adams and J. J. F. Fournier, \emph{Sobolev Spaces}, second edition, Academic Press, New York, 2003.

\bibitem{AHABW20}
F. Aldakheel, B. Hudobivnik, E. Artioli, L. Beir\~ao da Veiga, and P. Wriggers, Curvilinear virtual elements
for contact mechanics, \emph{Comput.\ Methods Appl.\ Mech.\ Engrg.} {\bf 372} (2020), paper no.\ 113394.

\bibitem{AM25}
W. Akram and M. T. Mohan, Mixed finite element method for a hemivariational inequality of s stationary 
convective Brinkman-Forchheimer Extended Darcy equations, {\tt arXiv}:2508.02797.

\bibitem{Ant83}  
S. S. Antman, The influence of elasticity in analysis: modern developments, 
\emph{Bulletin of the American Mathematical Society} {\bf 9}\,(3) (1983), 267--291.

\bibitem{AH2009}
K. Atkinson and W. Han, \emph{Theoretical Numerical Analysis: A Functional Analysis Framework},
third edition, Springer, New York, 2009.

\bibitem{BBH17}
M. Barboteu, K. Bartosz, and W. Han, Numerical analysis of an evolutionary variational-hemivariational  
inequality with application in contact mechanics, \emph{Comput.\ Methods Appl.\ Mech.\ Engrg.}
{\bf 318} (2017), 882--897.

\bibitem{BBHJ15}
M. Barboteu, K. Bartosz, W. Han, and T. Janiczko, Numerical analysis of a hyperbolic hemivariational 
inequality arising in dynamic contact, \emph{SIAM Journal on Numerical Analysis} {\bf 53} (2015), 527--550.

\bibitem{BHM19}
M. Barboteu, W. Han, and S. Mig\'{o}rski, On numerical approximation of a variational-hemivariational 
inequality modeling contact problems for locking materials, \emph{Comput.\ Math.\ Appl.} 
{\bf 77} (2019), 2894--2905.

\bibitem{Ba20}
K. Bartosz, Numerical analysis of a nonmonotone dynamic contact problem of a non-clamped piezoelectric 
viscoelastic body, \emph{Evol.\ Equ.\ Control Theory} {\bf 9} (2020), 961--980.

\bibitem{beirao2013basic}
L. Beir\~ao da Veiga, F. Brezzi, A. Cangiani, G. Manzini, L. D. Marini, and A. Russo, Basic
principles of virtual element methods, \emph{Math.\ Models Methods Appl.\ Sci.} {\bf 23} (2013), 119--214.

\bibitem{beirao2013virtual}
L. Beir\~ao da Veiga, F. Brezzi, and L. D. Marini, Virtual elements for linear elasticity problems,
\emph{SIAM J. Numer.\ Anal.} {\bf 51} (2013), 794--812.


\bibitem{BS2008}
S. C. Brenner and L. R. Scott, \emph{The Mathematical Theory of Finite Element Methods}, third edition,
Springer-Verlag, New York, 2008.

\bibitem{Br72}
H. Br\'ezis, Probl\`emes unilat\'eraux, {\em J. Math.\ Pures et Appl.} {\bf 51} (1972), 1--168.

\bibitem{Bre2011}
H. Br\'ezis, \emph{Functional Analysis, Sobolev Spaces and Partial Differential Equations}, Springer, New York, 2011.

\bibitem{CL2021}
S. Carl and V. K. Le, \emph{Multi-Valued Variational Inequalities and Inclusions}, Springer, New York, 2021.

\bibitem{CLM2007}
S. Carl, V. K. Le, and D. Motreanu, \emph{Nonsmooth Variational Problems and Their Inequalities:
Comparison Principles and Applications}, Springer, New York, 2007.


\bibitem{CHY23}
G. Caselli, M. Hensel, and I. Yousept, Quasilinear variational inequalities in ferromagnetic shielding:
well-posedness, regularity, and optimal control, \emph{SIAM J. Control Optim.} {\bf 61} (2023), 2043--2068.

\bibitem{CMMY23}
J. Cen, S. Mig\'{o}rski, C. Min, and J.-C. Yao,  Hemivariational inequality for contaminant 
reaction–diffusion model of recovered fracturing fluid in the wellbore of shale gas reservoir,
\emph{Commun.\ Nonlinear Sci.\ Numer.\ Simul.} {\bf 118} (2023), paper no.\ 107020.

\bibitem{CSWL23}
X.-L. Cheng, X. Shen, X. Wang, and K. Liang, A deep neural network-based method for solving obstacle problems,
\emph{Nonlinear Anal.\ Real World Appl.} {\bf 72} (2023), paper no.\ 103864.

\bibitem{CHR2023}
F. Chouly, P. Hild, and Y. Renard, \emph{Finite Element Approximation of Contact and Friction in Elasticity},
Birkh\"{a}user/Springer, Cham, 2023.

\bibitem{Ci1978}
P. G. Ciarlet, \emph{The Finite Element Method for Elliptic Problems}, North Holland, Amsterdam, 1978.

\bibitem{CHKW22}
M. Cihan, B. Hudobivnik, J. Korelc, and P. Wriggers, A virtual element method for 3D contact problems with 
non-conforming meshes, \emph{Comput.\ Methods Appl.\ Mech.\ Engrg.} {\bf 402} (2022), paper no.\ 115385.

\bibitem{Cl75}
F. H. Clarke, Generalized gradients and applications, \emph{Trans.\ Am.\ Math.\ Soc.} {\bf 205} (1975), 247--262.

\bibitem{Cl1983}
F. H. Clarke, \emph{Optimization and Nonsmooth Analysis}, Wiley, New York, 1983.

\bibitem{Cl2013}
F. H. Clarke, \emph{Functional Analysis, Calculus of Variations and Optimal Control}, Springer, London, 2013.

\bibitem{CLSW1998}
F. H. Clarke, Y.S. Ledyaev, R. J. Stern, and P. R. Wolenski, \emph{Nonsmooth Analysis and Control Theory},
Springer, New York, 1998.


\bibitem{DMP2}
Z. Denkowski, S. Mig\'orski and N.S. Papageorgiou, \emph{An Introduction to Nonlinear Analysis: Applications},
Kluwer Academic/Plenum Publishers, Boston, Dordrecht, London, New York, 2003.

\bibitem{DL1976}
G. Duvaut and J.-L. Lions, \emph{Inequalities in Mechanics and Physics}, Springer-Verlag, Berlin, 1976.

\bibitem{EJK2005}
C. Eck, J. Jaru\v sek, and M. Krbec, \emph{Unilateral Contact Problems: Variational Methods and Existence Theorems}, 
Pure and Applied Mathematics {\bf 270}, Chapman/CRC Press, New York, 2005.

\bibitem{ET1976}
I. Ekeland and R. Temam, \emph{Convex Analysis and Variational Problems}, North-Holland, Amsterdam, 1976.

\bibitem{Evans2010}
L. C. Evans, \emph{Partial Differential Equations}, second edition, American Mathematical Society, 2010.

\bibitem{Fa74}
R.S. Falk, Error estimates for the approximation of a class of variational inequalities,
\emph{Mathematics of Computation} {\bf 28} (1974), 963--971.

\bibitem{FCHCD20} C. Fang, K. Czuprynski, W. Han, X.L. Cheng, and X. Dai,
Finite element method for a stationary Stokes hemivariational inequality with slip boundary
condition, \emph{IMA J. Numer.\ Anal.} {\bf 40} (2020), 2696--2716.

\bibitem{FHH19}
F. Feng, W. Han, and J. Huang, Virtual element method for elliptic hemivariational inequalities,
\emph{Journal of Scientific Computing} {\bf 81} (2019), 2388--2412.

\bibitem{FHH21a}
F. Feng, W. Han, and J. Huang, Virtual element method for elliptic hemivariational inequalities with
a convex constraint, \emph{Numerical Mathematics: Theory, Methods and Applications} {\bf 14} (2021), 589--612.

\bibitem{FHH21b}
F. Feng, W. Han, and J. Huang, The virtual element method for an obstacle problem of a Kirchhoff plate,
\emph{Commun.\ Nonlinear Sci.\ Numer.\ Simul.}{\bf 103} (2021), paper no.\ 106008.

\bibitem{FHH22}
F. Feng, W. Han, and J. Huang, A nonconforming virtual element method for a fourth-order hemivariational inequality
in Kirchhoff plate problem, \emph{Journal of Scientific Computing} {\bf 90} (2022), paper no.\ 89.

\bibitem{Fi64}
G. Fichera, Problemi elastostatici con vincoli unilaterali: il problema di Signorini
con ambigue condizioni al contorno, \emph{Atti Accad.\ Naz.\ Lincei, Mem.,
Cl.\ Sci.\ Fis.\ Mat.\ Nat., Sez.\ I, VIII.\ Ser.} {\bf 7} (1964), 91--140.

\bibitem{Fu93}
H. Fujita, \emph{Flow Problems with Unilateral Boundary Conditions}, College de France, Lecons, 1993.

\bibitem{Fu94}
H. Fujita, A mathematical analysis of motions of viscous incompressible fluid under leak or slip boundary
conditions, \emph{RIMS K\^{o}ky\^{u}roku} {\bf 888} (1994), 199--216.

\bibitem{Gl1984}
R. Glowinski, \emph{Numerical Methods for Nonlinear Variational Problems}, Springer-Verlag, New York, 1984.

\bibitem{GLT1981}
R. Glowinski, J.-L. Lions, and R. Tr\'{e}moli\`{e}res, \emph{Numerical Analysis of Variational Inequalities},
North-Holland, Amsterdam, 1981.

\bibitem{GM2003}
D. Goeleven and D. Motreanu, \emph{Variational and Hemivariational Inequalities: Theory, Methods and
Applications. Vol.\ I. Unilateral Analysis and Unilateral Mechanics}, Kluwer Academic Publishers, Boston, MA, 2003.

\bibitem{GMDR2003}
D. Goeleven, D. Motreanu, Y. Dumont, and M. Rochdi, \emph{Variational and Hemivariational Inequalities: Theory,
Methods and Applications. Vol.\ II. Unilateral Problems}, Kluwer Academic Publishers, Boston, MA, 2003.

\bibitem{GJKR2022}
J. Gwinner, B. Jadamba, A. A. Khan, and F. Raciti, \emph{Uncertainty Quantification in Variational 
Inequalities: Theory, Numerics, and Applications}, CRC Press, Boca Raton, Florida, 2022.

\bibitem{HH19}
D. Han and W. Han, Numerical analysis of an evolutionary variational-hemivariational inequality with 
application to a dynamic contact problem, \emph{Journal of Computational and Applied Mathematics} 
{\bf 358} (2019), 163--178.

\bibitem{Han18}
W. Han, Numerical analysis of stationary variational-hemivariational inequalities with applications
in contact mechanics, \emph{Mathematics and Mechanics of Solids} {\bf 23} (2018), 279--293.

\bibitem{Han20}
W. Han, Minimization principles for elliptic hemivariational inequalities,
\emph{Nonlinear Analysis: Real World Applications} {\bf 54} (2020), paper no.\ 103114.

\bibitem{Han21}
W. Han, A revisit of elliptic variational-hemivariational inequalities,
\emph{Numerical Functional Analysis and Optimization} {\bf 42} (2021), 371--395.

\bibitem{Han2024}
W. Han, \emph{An Introduction to Theory and Applications of Stationary Variational-Hemivariational Inequalities},
Springer, New York, 2024.

\bibitem{HCJ21}
W. Han, K. Czuprynski, and F. Jing, Mixed finite element method for a hemivariational inequality of 
stationary Navier-Stokes equations, \emph{Journal of Scientific Computing} {\bf 89} (2021), paper no.\  8. 
        
\bibitem{HFWH25}
W. Han, F. Feng, F. Wang, and J. Huang, Numerical analysis of hemivariational inequalities with 
applications in contact mechanics, \emph{Advances in Applied Mechanics} {\bf 60} (2025), 113--178.

\bibitem{HJY23}
W. Han, F. Jing, and Y. Yao, Stabilized mixed finite element methods for a Navier--Stokes
hemivariational inequality, \emph{BIT Numerical Mathematics} {\bf 63} (2023), paper no.\  46.

\bibitem{HLW23}
W. Han, M. Ling, and F. Wang, Numerical solution of an $H({\rm curl})$-elliptic hemivariational 
inequality, \emph{IMA J. Numer.\ Anal.} {\bf 43} (2023), 976--1000.
        
\bibitem{HM22a}
W. Han and A. Matei, Minimax principles for elliptic mixed hemivariational-variational inequalities, 
\emph{Nonlinear Analysis: Real World Applications} {\bf 64} (2022), paper no.\  103448.

\bibitem{HM22b}
W. Han and A. Matei, Well-posedness of a general class of elliptic mixed hemivariational-variational
inequalities, \emph{Nonlinear Analysis: Real World Applications} {\bf 66} (2022), 103553.

\bibitem{HMS14}
W. Han, S. Mig\'orski and M. Sofonea, A class of variational-hemivariational inequalities with applications
to frictional contact problems, \emph{SIAM Journal on Mathematical Analysis} {\bf 46} (2014), 3891--3912.

\bibitem{HN23}
W. Han and M. Nashed, On variational-hemivariational inequalities in Banach spaces, \emph{Communications 
in Nonlinear Science and Numerical Simulation} {\bf 124} (2023), paper no,\ 107309.

\bibitem{HQM24}
W. Han,  H. Qiu, and L. Mei,   On a Stokes hemivariational inequality for incompressible fluid
flows with damping, \emph{Nonlinear Analysis: Real World Applications} {\bf 79} (2024), paper no.\ 104131.
        
\bibitem{HR2013}
W. Han and B.D. Reddy, \emph{Plasticity: Mathematical Theory and Numerical Analysis}, second edition, 
Springer-Verlag, 2013.

\bibitem{HS2002}
W. Han and M. Sofonea, \emph{Quasistatic Contact Problems in Viscoelasticity and Viscoplasticity},
American Mathematical Society and International Press, 2002.

\bibitem{HS19AN}
W. Han and M. Sofonea, Numerical analysis of hemivariational inequalities in contact mechanics,
\emph{Acta Numerica} {\bf 28} (2019), 175--286.

\bibitem{HSB17}
W. Han, M. Sofonea, and M. Barboteu,  Numerical analysis of elliptic hemivariational inequalities,
\emph{SIAM J. Numer.\ Anal.} {\bf 55} (2017), 640--663.

\bibitem{HSD18}
W. Han, M. Sofonea, and D. Danan, Numerical analysis of stationary variational-hemivariational
inequalities, \emph{Numer.\ Math.} {\bf 139} (2018), 563--592.

\bibitem{HZ19}
W. Han and S. Zeng, On convergence of numerical methods for variational-hemivariational inequalities
under minimal solution regularity, \emph{Applied Mathematics Letters} {\bf 93} (2019), 105--110.

\bibitem{HS66}
P. Hartman and G. Stampacchia, On some nonlinear elliptic differential functional equations,
\emph{Acta Math.} {\bf 15} (1966), 271--310.

\bibitem{HHN96}
J. Haslinger, I. Hlav\'{a}\v{c}ek,  J. Ne\v{c}as, Numerical methods for unilateral problems in solid mechanics,
in \emph{Handbook of Numerical Analysis, Vol.\ IV},  P.G. Ciarlet and J.L. Lions, eds.,
North-Holland, Amsterdam, 1996, 313--485.

\bibitem{HMP1999}
J. Haslinger, M. Miettinen and P. D. Panagiotopoulos, \emph{Finite Element Method for Hemivariational Inequalities:
Theory, Methods and Applications}, Kluwer Academic Publishers, Boston, Dordrecht, London, 1999.

\bibitem{HHNL1988} I. Hlav\'{a}\v{c}ek, J. Haslinger, J. Ne\v{c}as, and J. Lov\'{\i}\v{s}ek,
\emph{Solution of Variational Inequalities in Mechanics}, Springer-Verlag, New York, 1988.

\bibitem{HWW22}
J. Huang, C. Wang, and H. Wang, A deep learning method for elliptic hemivariational inequalities,
\emph{East Asian J. Appl.\ Math.} {\bf 12} (2022), 487--502.

\bibitem{JA2023}
A. Jayswal and T. Antczak (eds.), \emph{Continuous Optimization and Variational Inequalities},
CRC Press, Boca Raton, 2023.

\bibitem{JO21}
M. Jureczka and A. Ochal, A nonsmooth optimization approach for hemivariational inequalities with 
applications to contact mechanics, \emph{Appl.\ Math.\ Optim.} {\bf 83} (2021), 1465--1485.

\bibitem{JOB23}
M. Jureczka, A. Ochal, and P. Bartman, A nonsmooth optimization approach for time-dependent
hemivariational inequalities, \emph{Nonlinear Anal.\ Real World Appl.} {\bf 73} (2023), paper no.\ 103871.

\bibitem{KO1988} N. Kikuchi and J.T. Oden, \emph{Contact Problems in Elasticity:
A Study of Variational Inequalities and Finite Element Methods}, SIAM, Philadelphia, 1988.

\bibitem{LHZ22} M. Ling, W. Han, and S. Zeng, A pressure projection stabilized mixed finite element method for 
a Stokes hemivariational inequality, \emph{Journal of Scientific Computing} {\bf 92} (2022), paper no.\ 13.
    
\bibitem{LWH20}
M. Ling, F. Wang, and W. Han, The nonconforming virtual element method for a stationary Stokes hemivariational
inequality with slip boundary condition, \emph{Journal of Scientific Computing} {\bf 85} (2020), paper no.\ 56.

\bibitem{LXH24}
M. Ling, W. Xiao, and W. Han, Numerical analysis of a history-dependent mixed hemivariational-variational 
inequality in contact problems, \emph{Comput.\ Math.\ Appl.} {\bf 166} (2024), 65--76.

\bibitem{LS67}
J.-L. Lions and G. Stampacchia, Variational inequalities, \emph{Comm.\ Pure Appl.\ Math.} {\bf 20} (1967), 493--519.

\bibitem{MCHD24}
S. Mig\'orski, Y. Chao, J. He, and S. Dudek, Analysis of quasi-variational–hemivariational inequalities 
with applications to Bingham-type fluids, \emph{Commun.\ Nonlinear Sci.\ Numer.\ Simul.} {\bf 133} (2024), paper no.\ 107968.

\bibitem{MD22}
S. Mig\'orski and S. Dudek, A class of variational–hemivariational inequalities for Bingham type fluids,
\emph{Applied Mathematics \& Optimization} {\bf 85}, paper no.\ 16.

\bibitem{MOS2013}
S. Mig\'orski, A. Ochal, and M. Sofonea, \emph{Nonlinear Inclusions and Hemivariational
Inequalities. Models and Analysis of Contact Problems}, Advances in Mechanics
and Mathematics 26, Springer, New York, 2013.

\bibitem{MP1999}
D. Motreanu and P. D. Panagiotopoulos,  \emph{Minimax Theorems and Qualitative Properties of the
Solutions of Hemivariational Inequalities}, Kluwer Academic Publishers, Berlin, 1999.

\bibitem{NP1995}
Z. Naniewicz and P. D. Panagiotopoulos, \emph{Mathematical Theory of Hemivariational
Inequalities and Applications}, Dekker, New York, 1995.

\bibitem{NH1981} J. Ne\v cas and I. Hlava\v cek, \emph{Mathematical Theory of
Elastic and Elastoplastic Bodies: An Introduction}, Elsevier, Amsterdam, 1981.

\bibitem{OJB22} 
A. Ochal, M. Jureczka, and P. Bartman, A survey of numerical methods for hemivariational inequalities with
applications to contact mechanics, \emph{Commun.\ Nonlinear Sci.\ Numer.\ Simul.} {\bf 114} (2022), paper no.\ 106563.

\bibitem{Pa83}
P. D. Panagiotopoulos, Nonconvex energy functions, hemivariational inequalities
and substationary principles, \emph{Acta Mechanica} {\bf 42} (1983), 160--183.

\bibitem{Pa1993}
P. D. Panagiotopoulos, \emph{Hemivariational Inequalities, Applications in
Mechanics and Engineering}, Springer-Verlag, Berlin, 1993.

\bibitem{QLWW25}
J. Qiu, M. Ling, F. Wang, and B. Wu, Nonconforming finite element method for a 4th-order 
history-dependent hemivariational inequality, \emph{Commun.\ Nonlinear Sci.\ Numer.\ Simul.}
{\bf 145} (2025), paper no.\ 108750.

\bibitem{QWLZ23}
J. Qiu, F. Wang, M. Ling, and J. Zhao, The interior penalty virtual element method for the 
fourth-order elliptic hemivariational inequality, \emph{Commun.\ Nonlinear Sci.\ Numer.\ Simul.}
{\bf 127} (2023), paper no.\ 107547.

\bibitem{Si33}
A. Signorini, Sopra alcune questioni di elastostatica,
\emph{Atti della Societ\`a Italiana per il Progresso delle Scienze}, 1933.

\bibitem{SHS2006}
M. Sofonea, W. Han, and M. Shillor, \emph{Analysis and Approximation of Contact Problems with
Adhesion or Damage}, Chapman \& Hall/CRC, New York, 2006.

\bibitem{SM2009}
M. Sofonea and A. Matei, \emph{Variational Inequalities with Applications: A Study of Antiplane Frictional
Contact Problems}, Springer, 2009.

\bibitem{SM11}
M. Sofonea and A. Matei, History-dependent quasivariational inequalities arising in contact mechanics,
\emph{European Journal of Applied Mathematics} {\bf 22} (2011), 471--491.

\bibitem{SM2012}
M. Sofonea and A. Matei, \emph{Mathematical Models in Contact Mechanics}, Cambridge University Press, Cambridge, 2012.

\bibitem{SM16}
M. Sofonea and S. Mig\'orski, A class of history-dependent variational-hemivariational inequalities,
\emph{Nonlinear Differ.\ Equ.\ Appl.} {\bf 23} (216), paper no.\ 38.

\bibitem{SM2025}
M. Sofonea and S. Mig\'orski, \emph{Variational-Hemivariational Inequalities with Applications},
second edition, CRC Press, Boca Raton, FL, 2025.

\bibitem{St64}
G. Stampacchia,  Formes bilineaires coercitives sur les ensembles convexes,
\emph{C. R. Acad.\ Sci.} {\bf 258} (1964), 4413--4416.

\bibitem{SCX24}
X. Sun, X.-L. Cheng, and H. Xuan, Error estimates and numerical simulations of a thermoviscoelastic 
contact problem with damage and long memory, \emph{Commun.\ Nonlinear Sci.\ Numer.\ Simul.}
{\bf 137} (2024), paper no.\ 108165.

\bibitem{TC25}
X. Tan and T. Chen, A mixed finite element approach for a variational-hemivariational inequality of 
incompressible Bingham fluids, \emph{J. Sci.\ Comput.} {\bf 103} (2025), paper no.\ 36.

\bibitem{Te1985}
R. Temam, \emph{Mathematical Problems in Plasticity}, Gauthier-Villlars, Paris, 1985.

\bibitem{UL2011}
M. Ulbrich, \emph{Semismooth Newton Methods for Variational Inequalities and Constrained Optimization
Problems in Function Spaces}, SIAM, 2011.

\bibitem{WD24}
F. Wang and H. Dang, Randomized neural network methods for solving obstacle problems, 
\emph{Banach Center Publ.} {\bf 127} (2024), 261--276.

\bibitem{WQ20}
F. Wang and H. Qi, A discontinuous Galerkin method for an elliptic hemivariational
inequality for semipermeable media, \emph{Applied Mathematics Letters} {\bf 109} (2020), paper no.\ 106572. 

\bibitem{WSW23}
F. Wang, S. Shah, and B. Wu, Discontinuous Galerkin methods for hemivariational inequalities in
contact mechanics, \emph{J. Sci.\ Comp.} {\bf 95} (2023), paper no.\ 87.

\bibitem{WWH21}
F. Wang, B. Wu, and W. Han, The virtual element method for general elliptic hemivariational inequalities,
\emph{Journal of Computational and Applied Mathematics} {\bf 389} (2021), paper no.\ 113330.

\bibitem{WZ21}
F. Wang and J. Zhao, Conforming and nonconforming virtual element methods for a Kirchhoff 
plate contact problem, \emph{IMA J. Numer.\ Anal.} {\bf 41} (2021), 1496-–1521. 

\bibitem{WXHC20}
S. Wang, W. Xu, W. Han, and W. Chen, Numerical analysis of history-dependent variational-hemivariational 
inequalities, \emph{Sci.\ China Math.} {\bf 63} (2020), 2207–-2232.

\bibitem{WCH26}
W. Wang, X.-L. Cheng, and W. Han, Analysis and finite element solution of a Navier–Stokes hemivariational 
inequality for incompressible fluid flows with damping, \emph{Nonlinear Analysis: Real World Applications} 
{\bf 87} (2026), paper no.\ 104439.

\bibitem{Wr2006}
P. Wriggers, \emph{Computational Contact Mechanics}, second edition, Springer, Berlin, 2006.

\bibitem{wriggers2016virtual}
P. Wriggers, W. T. Rust, B. D. Reddy, A virtual element method for contact,
\emph{Comput.\ Mech.} {\bf 58} (2016), 1039--1050.

\bibitem{WWH22}
B. Wu, F. Wang, and W. Han, Virtual element method for a frictional contact problem with normal compliance,
\emph{Commun.\ Nonlinear Sci.\ Numer.\ Simul.} {\bf 107} (2022), paper no.\ 106125.

\bibitem{WWH24}
B. Wu, F. Wang, and W. Han, The virtual element method for a contact problem with wear and unilateral constraint,
\emph{Appl.\ Numer.\ Math.} {\bf 206} (2024), 29--47. 

\bibitem{XL23a}
W. Xiao and M. Ling, The virtual element method for general variational-hemivariational inequalities
with applications to contact mechanics, \emph{J. Comput.\ Appl.\ Math.} {\bf 428} (2023), paper no.\ 115152.

\bibitem{XL23b}
W. Xiao and M. Ling, A priori error estimate of virtual element method for a quasivariational-hemivariational
inequality, \emph{Commun.\ Nonlinear Sci.\ Numer.\ Simul.} {\bf 121} (2023), paper no.\ 107222.

\bibitem{XL23}
W. Xiao and M. Ling, Virtual element method for a history-dependent variational-hemivariational inequality 
in contact problems, \emph{J. Sci.\ Comput.} {\bf 96} (2023), paper no.\ 82.

\bibitem{XHHCW20}
W. Xu, Z. Huang, W. Han, W. Chen, and C. Wang, Numerical approximation of an electro-elastic frictional 
contact problem modeled by hemivariational inequality, \emph{Comput.\ Appl.\ Math.} {\bf 39} (2020), 265.

\bibitem{XHHCW19}
W. Xu, Z. Huang, W. Han, W. Chen, and C. Wang, Numerical analysis of history-dependent variational-hemivariational 
inequalities with applications in contact mechanics, \emph{J. Comput.\ Appl.\ Math.} {\bf 351} (2019), 364--377.

\bibitem{XCY25}
H. Xuan, X.-L. Cheng, and L. Yuan, Numerical studies of a class of thermoviscoelastic frictional contact 
problem described by fractional differential hemivariational inequalities, \emph{J. Sci.\ Comput.}
{\bf 103} (2025), paper no.\ 4.

\bibitem{Yo21}
I. Yousept, Maxwell quasi-variational inequalities in superconductivity,
\emph{ESAIM Math.\ Model.\ Numer.\ Anal.} {\bf 55} (2021), 1545--1568.

\bibitem{ZeIII}
E. Zeidler, \emph{Nonlinear Functional Analysis and its Applications. III: Variational Methods and Optimization}, 
Springer-Verlag, New York, 1985.

\end{thebibliography}
\end{document}